\documentclass[reqno,12pt]{amsart}

\newtheorem{theorem}{Theorem}[section]
\newtheorem{lemma}[theorem]{Lemma}
\newtheorem{corollary}[theorem]{Corollary}

\theoremstyle{definition}

\newtheorem{assumption}[theorem]{Assumption}

\theoremstyle{remark}
\newtheorem{remark}[theorem]{Remark}

\newcommand{\mysection}[1]{\section{#1}
\setcounter{equation}{0}}

\newcommand{\bR}{\mathbb R}
\newcommand{\bZ}{\mathbb Z}

\newcommand\cM{\mathcal{M}}
\newcommand\cK{\mathcal{K}}

\begin{document}
\title[On the rate of convergence] {On the rate of convergence
of finite-difference approximations for Bellman equations with
Lipschitz coefficients}

\author[N.V. Krylov]{N.V. Krylov}

\address
{127 Vincent Hall, University of Minnesota, Minneapolis, MN
55455, USA}

\email{krylov@math.umn.edu}



\begin{abstract}
We consider parabolic Bellman equations 
with Lipschitz coefficients. Error bounds 
of order $h^{1/2}$ for certain
types of finite-difference schemes are obtained. 

\end{abstract}

\maketitle

\mysection{Introduction}

Bellman equations arise in many areas of mathematics, say
in control theory, differential geometry,
and mathematical finance, to name a few.
These equations typically are fully nonlinear
second order degenerate elliptic or parabolic equations.
In the particular case of complete degeneration
they become Hamilton-Jacobi first-order equations.

Quite naturally, the problem of finding
numerical methods of approximating
solutions to Bellman equations arises.
First methods
dating back some thirty years ago
 were based on the fact that
the solutions are the value functions
in certain problems for controlled  diffusion processes,
that can be approximated by controlled Markov
chains. An account of the results obtained in this direction
can be found in  \cite{KD} and \cite{FS}. 

Another approach is based on the notion
of viscosity solution, which allows one to avoid
using probability theory. We refer to
\cite{BJ1} and \cite{BJ2} and the references
therein for discussion of what is achieved
in this direction. 

We will be dealing with degenerate {\em second-order\/}
equations. There is a very extensive literature
treating  Hamilton-Jacobi equations
and establishing the rate of convergence
of various numerical approximations.
The reader can find how much was done
for them in \cite{BC} and \cite{FG}.
In contrast, until quite recently there were
no results about the rate of convergence
of finite-difference approximations
for degenerate Bellman equations.
The first result appeared only in 1997 
for elliptic Bellman equations with
constant ``coefficients" (see \cite{Kr97})
and they were later extended to variable coefficients
and parabolic equations
in \cite{BJ1},
\cite{BJ2}, \cite{Kr99},  and \cite{Kr00}.
Surprisingly, as far as we know until
now these are {\em the only\/} published
result on the rate of convergence
of finite-difference approximations
even if Bellman equation becomes a {\em linear\/}
second order
  degenerate equation. One has to notice
however that there is vast literature
about other type of numerical approximations
for linear degenerate equations such as Galerkin or
finite element approximations 
(see, for instance, \cite{LT}). It is also worth noting
that under variety of conditions
 the first {\em sharp\/} estimates
for finite-difference approximations
in linear one-dimensional degenerate
case are proved 
in \cite{Zh}.

Our approach is based on two ideas from
\cite{Kr97}, \cite{Kr99}, and \cite{Kr00}
that the original equation and its finite-difference
approximation should play symmetric roles
and that one can ``shake the coefficients"
of the equation in order to be able to
mollify under the sign of nonlinear operator.
While shaking the coefficients of the approximate equation
we encounter a major problem of estimating
how much the solution of the shaken equation
differs from the original one. Solving
this problem amounts to estimating the
Lipschitz constant of the approximate solution.
We prove this estimate on the 
basis of Theorem \ref{theorem 8.19.1}
and consider this theorem as the most
important technical result of the present paper.
Theorem \ref{theorem 8.19.1} is new even
if the equation is linear although in that
case one can give a much simpler proof,
which we intend to do in a subsequent
joint article with Hogjie Dong.

Our main result says that for parabolic equations
in a special form with $C^{1/2,1}$
coefficients
the rate of convergence is not less than
$\tau^{1/4}+h^{1/2}$, where $\tau$ and $h$
are the time and space steps, respectively.
Simple examples show that 
under our conditions the estimate is sharp
even for the  case of linear first order equations.
For the elliptic case the rate becomes
$h^{1/2}$, which under comparable
conditions is slightly better than $h^{1/5}$
from \cite{BJ2}.

 The main emphasis of this paper is
on {\em constructing\/} 
finite-difference
approximations as good as possible for a given Bellman equation.
There is another part of the story
when one is interested in how more or less
arbitrary consistent finite-difference
type approximations converge to the true
solution. In this directions the known results are
somewhat weaker. We only know that
for $\tau=h^{2}$ there is an estimate
of order $h^{1/21}$, which sometimes becomes
$h^{1/3}$ (see \cite{Kr97}, \cite{Kr99}).

One particular degenerate Bellman equation is worth
mentioning separately. This equation arises as an obstacle
problem in PDEs or as an optimal stopping 
problem in stochastic control:
$$
\max(\Delta u-u,-u+g)=0,
$$
where $g$ is a given function. One usually
rewrites it in an equivalent form:
$$
\Delta u-u\leq0,\quad g\leq u,\quad
\Delta u-u=0\quad\text{on}\quad\{u>g\}.
$$

To conclude the introduction, we introduce some notation:
$\bR^d$ is a $d$-dimensional Euclidean space ;
$x=(x^1,x^2,\dots,x^d)$ is a typical point in $\bR^d$.
For any $l \in \bR^d$ and any
differentiable function $u$ on $\bR^d$, we denote
$$
D_lu=u_{x^i}l^i,\quad D^2_lu=u_{x^ix^j}l^il^j,
$$ etc. 
The symbols $D^{n}_{t}u$ stand for   the $n$th derivative
 in $t$  of 
$u=u(t,x)$, $t\in\bR$, $x\in\bR^{d}$,
and $D^{n}_{x}u$ for the collection of all $n$th
order derivatives of $u$ in $x$. We also use 
the notation 
$$
|u|_{0,Q}=\sup_{Q}|u|.
$$

Various
constants are denoted by $N$ in general and the expression
$N=N(\cdots)$ means that the given constant $N$ depends only on
the contents of the parentheses. We set
$$
a_{\pm}=a^{\pm}=(1/2)(|a|\pm a).
$$

Finally,
as usual the summation convention over repeated indices
is enforced.

{\bf Acknowledgment}. 
 The work  was partially supported by
NSF Grant DMS-0140405. The author is sincerely
grateful to Hongjie Dong who made valuable
comments on the first version of the article.

\mysection{Main results}
                                       \label{section 10.4.1}

Let $A$ be a separable metric space, constants 
$$
T \in(0,\infty),\quad
 K\in[1,\infty),\quad
\lambda\in[0,\infty),\quad
\text{integers}\quad  d,d_{1}\geq1.
$$

Suppose that we are given $\ell_{k}\in\bR^{d}$
and real-valued
$$\sigma^{\alpha}_{k}(t, x),\quad
b^{\alpha}_{k}(t, x),\quad
c^{\alpha} (t, x),\quad
f^{\alpha} (t, x),\quad
g(x)
$$ defined for
$k=\pm1,...,\pm d_{1}$, $(t,x)\in\bR\times \bR^{d}$, and
$\alpha\in A$ 
such that 
$$
\ell_{k}=-\ell_{-k},\quad\sigma^{\alpha}_{k}=\sigma^{\alpha}_{-k},
\quad b^{\alpha}_{k}\geq0,
\quad  c^{\alpha}\geq\lambda ,\quad|\ell_{k}|\leq K,
$$
\begin{assumption}
                                   \label{assumption 8.25.1}

 For $\psi 
= \sigma^{\alpha}_{k},b^{\alpha}_{k},c^{\alpha}
-\lambda,
f^{\alpha},g $, $k=1,...,d_{1}$,
$\alpha\in A$, for each $t\in[0,T]$  and $x,y\in\bR^{d}$ we have
$$
|\psi (t,x)|\leq K,\quad|\psi(t,x)-\psi(t,y)|\leq
K|x-y|.
$$
\end{assumption}

We also assume that $\sigma^{\alpha}_{k}(t, x)$, 
$b^{\alpha}_{k}(t, x)$, 
$c^{\alpha} (t, x) $, $f^{\alpha} (t, x)$
are Borel in $t$ and continuous in $\alpha$.

Introduce
$$
F(p_{k},q_{k},r,t,x)
=\sup_{\alpha\in A}[a^{\alpha}_{k}(t,x)p_{k}
+b^{\alpha}_{k}(t,x)q_{k}-c^{\alpha}(t,x)r+f^{\alpha}(t,x)]
$$
with the summation in $k$ performed before the supremum is taken.

Under the above assumptions there is a probabilistic
solution $v$ of the Bellman equation
\begin{equation}
                                                    \label{8.19.3}
\frac{\partial}{\partial t} u(t,x)+
F(D^{2}_{\ell_{k}}u(t,x),D_{\ell_{k}}u(t,x),
u(t,x),t,x)
=0 
\end{equation}
in 
$$
H_{T}:=[0,T)\times\bR^{d}
$$
 with terminal condition
\begin{equation}
                                                    \label{9.11.4}
u(T,x)=g(x),\quad x\in\bR^{d}.
\end{equation}
This solution is constructed by means of control theory.
The  reader unfamiliar with control theory
may consider $v$ as the unique
bounded viscosity solution
of the above problem   (see, for instance, \cite{FS},
\cite{KD}).

 For $h,\tau>0$, $l\in\bR^{d}$, $(t,x)\in[0,T)\times\bR^{d}$
 introduce
$$
\delta_{h,l}u(t,x)=\frac{u(t,x+hl)-u(t,x)}{h},\quad
\Delta_{h,l}u(t,x)=-\delta_{h,l}\delta_{h,-l}u(t,x)
$$
$$
=\frac{\delta_{h,l}
+\delta_{h,-l}}{h}u(t,x)
=\frac{u(t,x+hl)-2u(t,x)+u(t,x-hl)}{h^{2}},
$$
$$
\delta_{\tau}u(t,x)=\frac{u(t+\tau,x )-u(t,x)}{\tau},
$$
$$
\delta^{T}_{\tau}u(t,x)=\frac{u(t+\tau_{T}(t)
,x )-u(t,x)}{\tau },\quad\tau_{T}(t)=\tau\wedge(T-t).
$$
Just in case, notice that in the denominator
of $\delta^{T}_{\tau}u$ we write $\tau$ and not $\tau_{T}(t)$.
This will be important in the proof of Lemma \ref{lemma 9.13.3}.
Also note that
$$
t+\tau_{T}(t)=(t+\tau)\wedge T.
$$

Set
$$
\delta_{0,l}u=0,\quad
a^{\alpha}_{k}=(1/2)(\sigma^{\alpha}_{k})^{2}.
$$
In $H_{T}$
consider the following equation
with respect to a function $u$ given in $\bar{H}_{T}$
\begin{equation}
                                                 \label{8.19.1}
\delta^{T}_{\tau}u(t,x)+F(\Delta_{h,\ell_{k}}u(t,x),
\delta_{h,\ell_{k}}u(t,x),u(t,x),t, x)=0 
\end{equation}
with terminal condition \eqref{9.11.4}. 

Equation (\ref{8.19.1}) is an implicit finite-difference
approximation for the Bellman equation \eqref{8.19.3}.
Existence of a unique bounded solution
of problem \eqref{8.19.1}-\eqref{9.11.4}, which we denote
by $v_{\tau,h}$, is a standard fact proved by
successive approximations  in
Lemma \ref{lemma 9.7.1} (also see the comments before
that lemma).

Here are our main results.

\begin{theorem}
                                         \label{theorem 8.26.1}
In addition to the above assumptions
suppose that 

(H) for $\psi 
= \sigma^{\alpha}_{k},b^{\alpha}_{k},c^{\alpha}
-\lambda,
f^{\alpha} $, $k=1,...,d_{1}$,
$\alpha\in A$, for each 
  $x \in\bR^{d}$ and $t,s\in\bR$   we have
$$
 |\psi(t,x)-\psi(s,x)|\leq
K|t-s|^{1/2};
$$

Then there exists a constant $N_{1}$ depending only on
$d,d_{1},T$, and $K$ (but not $h$ or $\tau$)
 such that 
\begin{equation}
                                                  \label{9.18.2}
|v-v_{\tau,h}|\leq N_{1}(\tau^{1/4}+h^{1/2})
\end{equation}
in $H_{T}$. In addition, 
there exists a constant $N_{2}$ depending only on
$ d_{1}$  and $K$, such that
if $\lambda\geq N_{2}$, then $N_{1}$ is independent
of $T$.
\end{theorem}

Introduce 
$$
L^{\alpha}u=
a^{\alpha}_{k} D^{2}_{ \ell_{k}}u 
+b^{\alpha}_{k} 
D_{ \ell_{k}}u -
c^{\alpha} u ,\quad
L^{\alpha}_{h}u=
 a^{\alpha}_{k}\Delta_{h,\ell_{k}}u
+b^{\alpha}_{k}\delta_{h,\ell_{k}}u
-c^{\alpha}u.
$$

\begin{theorem}
                                         \label{theorem 8.27.1}
Under the assumptions before Theorem \ref{theorem 8.26.1}
suppose that $\sigma,b,c,f$ are independent of $t$
and $\lambda\geq N_{2}$, where $N_{2}$
is taken from Theorem \ref{theorem 8.26.1}.
Let $\tilde{v}(x)$ be a probabilistic or
the unique bounded viscosity solution of
$$
\sup_{\alpha\in A}
[L^{\alpha}u+f^{\alpha} ]=0
$$
in $\bR^{d}$. Let $\tilde{v}_{h}$ be the unique
bounded solution of
\begin{equation}
                                                 \label{8.27.2}
\sup_{\alpha\in A}
[L^{\alpha}_{h}u+f^{\alpha} ]=0
\end{equation}
in $\bR^{d}$. Then 
$$
|\tilde{v}-\tilde{v}_{h}|
\leq Nh^{1/2} 
$$
in $\bR^{d}$,  where $N$ depends only 
 on
$d,d_{1}$, and $K$.
\end{theorem}

The following result  about semidiscretization
allows one to use  
 approximations of the time derivative
different from the one in \eqref{8.19.1},
in particular, explicit schemes could be used.
\begin{theorem}
                                          \label{theorem 9.19.1} 
Under the assumptions of
Theorem \ref{theorem 8.26.1} there exists a unique
bounded solution $v_{h}(t,x)$ of 
\begin{equation}
                                                 \label{9.19.1}
\frac{\partial}{\partial t}u(t,x)+F(\Delta_{h,\ell_{k}}u(t,x),
\delta_{h,\ell_{k}}u(t,x),u(t,x),t, x)=0
\end{equation}
in $H_{T}$ with terminal condition \eqref{9.11.4}.
Furthermore, there exists a constant $N_{1}$ depending
only on $K$, $T$, $d$, and $d_{1}$ such that
$$
|v -v_{h} |\leq N_{1}h^{1/2} 
$$
in $H_{T}$. Finally,
there is a constant $N_{2}$ depending only
on $K$ and $d_{1}$ such that if $\lambda\geq N_{2}$,
then $N_{1}$ is independent of $T$.

\end{theorem}

We prove the above results in Section \ref{section 10.2.1},
after proving some auxiliary statements in
Sections \ref{section 9.10.1} and \ref{section 10.2.3}.
Then come  the main estimate of the Lipschitz constant in $x$
in Section \ref{section 10.2.4}
and finally  the H\"older 1/2 continuity in $t$
in Section \ref{section 10.2.3}.

\begin{remark}
In a subsequent article we will show
that assumption (H) is not needed
in Theorem \ref{theorem 9.19.1}.
Few other possible extensions of the above results
are discussed in Section~\ref{section 9.19.1}.
\end{remark}

\begin{remark}
                                      \label{remark 4.11.1}
One may think that considering the 
operators $L^{\alpha}$
 written
in the form $a_{k}^\alpha D^{2}_{\ell_{k}}
+b_{k}^\alpha D_{\ell_{k}}+c^{\alpha}$
is a severe restriction.
However it is easy to see (cf.~\cite{DK}) that if we fix a 
finite subset $B\subset\bZ^{d}$, such that
$\text{Span}\,B=\bR^{d}$,
 and if an operator $Lu=a^{ij}u_{x^{i}x^{j}}
+b^{i}u_{x^{i}}$ admits a  finite-difference approximation
$$
L_{h}u(x)=\sum_{y\in B}p_{h}(y)u(x+hy)\to Lu(x)
\quad\forall u\in C^{2}
$$
and $L_{h}$ are 
  monotone, then automatically
$$
L=\sum_{\stackrel{\scriptstyle l\in B}{l\ne0}}a_{l}
D^{2}_{l }+\sum_{\stackrel{\scriptstyle l\in
B}{l\ne0}}b_{l}D_{l}
$$
for some $a_{l}\geq0$ and $b_{l}\in\bR$.

There is also a very substantial
 advantage of using this particular 
form of $L^{\alpha}$ because 
 for any smooth function $\eta$ by Taylor's
formula we have
$$
D^{ 2}_l\eta(y)
=\Delta^{2}_{h,l}\eta(y)
-\frac{1}{6h^{2} }\int_{-h}^{h}D^{ 4}_l\eta(y+sl)
(h-|s|)^{3}\,ds
$$
and the second term on the right
has order $h^{2}$. By considering similarly
first order terms we see that
for any four times continuously differentiable function
$\eta$
\begin{equation}
                                       \label{9.9.2}
|L^{\alpha}\eta(x)-L^{\alpha}_{h}\eta(x)|
\leq N^{*}(h^{2}\sup_{B_{K}(x)} |D^{4}_{x}\eta|
+h \sup_{B_{K}(x)} |D^{2}_{x}\eta|),
\end{equation}
where $B_{K}(x)$ is the ball of radius $K$ centered at $x$
and $N^{*}$ depends only on $K$ and $d_{1}$.

\end{remark}

 \mysection{Solvability and
comparison principle for finite-difference
equations}
                                        \label{section 9.10.1}

Problem \eqref{8.19.1}-\eqref{9.11.4}
is, actually, a collection
of disjoint problems given on each mesh
associated with points $(t_{0},x_{0})\in[0,T)\times\bR^{d}$:
$$
\{\big((t_{0}+j\tau)\wedge T,
x_{0}+h(i_{1}\ell_{1}+...+i_{d_{1}}\ell_{d_{1}})\big)
 :
$$
\begin{equation}
                                             \label{9.9.5}
j=0,1,...,i_{k}=0,\pm1,...,
k=1,...,d_{1}\}.
\end{equation}
Indeed, problem 
\eqref{8.19.1}-\eqref{9.11.4} on such a mesh has perfect
sense even if $u$ is defined only on it.
In the future we will see that it is extremely convenient
to consider this collection of problems simultaneously.
However, while obtaining certain estimates
it is more convenient to work in a more traditional setting
with each particular mesh separately. In this way
even the  results look more general
and the continuity hypothesis in $t$ on the coefficients often
becomes just superfluous. It is also worth noting that
we do not assume that $\{\ell_{k}\}$ generates
$\bR^{d}$ so that the meshes (\ref{9.9.5}) may be
meshes on hyperplanes.

For fixed
$\tau,h>0$   introduce
$$
\bar\cM_{T}=\{(t,x)\in[0,T]
\times\bR^{d}:t=(j\tau)\wedge T,x=h(i_{1}\ell_{1}+
...+i_{d_{1}}\ell_{d_{1}}),
$$
$$
j=0,1,..., i_{k}=0,\pm1,....,k=1,...,d_{1}\}.
$$
Of course,   results obtained for equations
on subsets of $\bar\cM_{T}$ automatically translate
into the corresponding results for all other
meshes like~\eqref{9.9.5}.

Take a nonempty set 
$$
Q\subset\cM_{T} :=\bar\cM_{T}\cap([0,T)\times\bR^{d}).
$$
We start with a solvability result.
\begin{lemma}
                                             \label{lemma 9.7.1}

  Let $g(t,x)$ be a bounded function
on $\bar\cM_{T}$.  Then
there is a unique bounded function $u$
defined on $\bar\cM_{T}$
such that
equation \eqref{8.19.1} holds in $Q$ and
 $u=g$ on $\bar\cM_{T}\setminus Q$.

\end{lemma}

Proof.  Take a  constant $\gamma\in(0,1)$ 
  and define a function $\xi(t)
=\xi(t,x)$ on $\bar\cM_{T }$ recursively by
\begin{equation}
                                                \label{9.31.1}
\xi(T)=1,\quad\xi(t)=\gamma^{-1}\xi(t+\tau_{T}(t) )\quad
\text{for}\quad t<T.
\end{equation}
Notice that for any function $v$
\begin{equation}
                                                 \label{10.2.1}
\delta_{\tau}^{T}(\xi v) =\gamma\xi  \delta^{T}_{\tau}
v-\nu \xi v  ,\quad
\nu=\frac{1-\gamma}{\tau }.
\end{equation}
 Obviously the function $u$ we are looking for
is to satisfy
$$
u=\xi v,
$$
\begin{equation}
                                                 \label{9.8.1}
v(t,x)=\xi^{-1}(t )g (t,x)I_{\bar\cM_{T}\setminus Q}(t,x)
+I_{Q}(t,x)G[v] (t,x),
\end{equation}
where
$$
G[v](t,x) := v(t,x)+\varepsilon\xi^{-1}(t )[
\delta^{T}_{\tau}u(t,x)
$$
$$
+F(\Delta_{h,\ell_{k}}u(t,x),
\delta_{h,\ell_{k}}u(t,x),u(t,x),t, x)]
$$
and $\varepsilon$ is any number. Observe that
for $\varepsilon>0$
$$
G[v] (t,x)=\sup_{\alpha\in A}
[p_{\tau}v(t+\tau_{T}(t),x)
+p^{\alpha}_{k}(t,x) v(t,x+h\ell_{k}) 
$$
\begin{equation}
                                                 \label{9.31.2} 
+p^{\alpha}(t,x)
v(t,x)+\varepsilon\xi^{-1}(t+\tau_{T}(t))f^{\alpha}(t,x)],
\end{equation}
where
$$
p_{\tau}=\varepsilon\gamma\tau^{-1},\quad
p^{\alpha}_{k} =2
\varepsilon h^{-2}a^{\alpha}_{k} 
+\varepsilon h^{-1}b^{\alpha}_{k} ,
$$
$$
p^{\alpha} :=
 1-p_{\tau}- 
\sum_{k}p^{\alpha}_{k} 
 -\varepsilon\nu-
\varepsilon c^{\alpha} .
$$
We choose $\varepsilon$ and $\gamma$ so that
$$
 p^{\alpha}_{k},p^{\alpha}\geq0,\quad
0\leq\sum_{k}p^{\alpha}_{k}+p^{\alpha}+p_{\tau}=
1 -\varepsilon\nu-
\varepsilon c^{\alpha}\leq\delta<1,
$$
where $\delta$ is a constant.

Then we use the fact that the difference of sups is less
than the sup of differences and easily conclude that 
for  any functions $v$ and $w$ we have
$$
|G[v](t,x)-G[w](t,x)|\leq\delta\sup_{\bar\cM_{T}}
|v-w|,
$$
so that the operator $G$ is a contraction in the space
of bounded functions on $ \bar\cM_{T}$. The application
of Banach's fixed point theorem
to equation \eqref{9.8.1} proves the lemma.

\begin{remark}
                                        \label{remark 9.17.1}
Sometimes dealing with functions
on $\bar\cM_{T}$  the fact that
 $T$ may not be a point of type $\tau,2\tau,...$
is quite inconvenient
just because then we should  take care of two cases: $t<T$
and $t=T$, separately. In addition, on few occasions
in the article we are not  
using any continuity hypotheses in $t$. Therefore,
we may move the points $(j\tau)\wedge T$ 
along the time axis preserving their order
in any way we like provided that we carry along with them
the values of the coefficients
and other functions involved.
In connection with this we introduce
$T'$ as the least point in the progression
$\tau,2\tau,...$, which is   $\geq T$
and notice that equation \eqref{8.19.1}
on $Q$ 
is rewritten as the following equation
on $Q$ relative to a function $\tilde{u}$
given on $\bar\cM_{T'}$:
$$
\delta_{\tau}\tilde{u}(t,x)+\sup_{\alpha\in A}
[L^{\alpha}_{h}(t,x)\tilde{u}(t,x)+f^{\alpha}(t,x)]=0
$$
where $\tilde{u}(t,x)=u(t,x)$ on $\cM_{T'}$ and
$\tilde{u}(T',x)=u(T,x)$. Observe that
$$
\delta_{\tau}\tilde{u}(t,x)=
\delta^{T'}_{\tau}\tilde{u}(t,x)=
\delta_{\tau}^{T}u(t,x)
$$
 on $\cM_{T'}$. Also note that the condition
  $u=g$ on $\bar\cM_{T}\setminus Q$ translates
into $\tilde{u}=\tilde{g}$ on $\bar\cM_{T'}\setminus Q$,
where $\tilde{g}(t,x)=g(t,x)$  on $\cM_{T'}$
and $\tilde{g}(T',x)=g(T,x)$.
\end{remark}

The following is a comparison result.

\begin{lemma}
                                             \label{lemma 9.9.2}
  Let  $u_{1},u_{2}$ be  functions on $\bar\cM_{T}$,
$f_{1}^{\alpha}(t,x),f_{2}^{\alpha}(t,x)$   functions
on $A\times\cM_{T}$ and $C$ a constant. Assume that in $Q$
$$
\sup_{\alpha}f_{2}^{\alpha}<\infty,\quad
f_{1}^{\alpha}\leq f_{2}^{\alpha},
$$
\begin{equation}
                                                 \label{9.8.2}
\delta^{T}_{\tau}u_{1}+ \sup_{\alpha\in A}[L^{\alpha}_{h}u_{1}+f_{1}^{\alpha}]
+C\geq \delta^{T}_{\tau}u_{2}+\sup_{\alpha\in A}[L^{\alpha}_{h}
u_{2}+f_{2}^{\alpha}].
\end{equation}
Finally, let $h\leq1$ and $u_{1}\leq u_{2}$ on 
$\bar\cM_{T}\setminus Q$ and
assume that $u_{i}e^{-\mu|x|}$ are bounded
on $\cM$, where $\mu\geq0$ is a constant.
We assert that there exists a  constant $\tau^{*}>0$,
depending only on $K$, $d_{1}$, and $\mu$,
such that if $\tau\in(0,\tau^{*})$ then on $\bar\cM_{T}$
\begin{equation}
                                                 \label{9.8.3}
u_{1}\leq u_{2}+  T' C_{+}.
\end{equation}

Furthermore, $\tau^{*}(K,d_{1},\mu)\to\infty$
as $\mu\downarrow0$ and
if $u_{1},u_{2}$ are bounded on $\bar\cM_{T}$,
so that $\mu=0$,
then \eqref{9.8.3} holds without any constraints on
$h$ and $\tau$.
\end{lemma}

Proof. Obviously, one can replace $f^{\alpha}_{1}$
with $f^{\alpha}_{2}$ preserving \eqref{9.8.2}.
Then, according   to Remark
\ref{remark 9.17.1}, we can pass from
  $T$ to $T'$ and thereby we may assume that
$T=T'$. We get from
\eqref{9.8.2} that
$$
\delta _{\tau}u+
\sup_{\alpha\in A}L^{\alpha}_{h}u+C\geq0
$$
on $Q$, where  $u=u_{1}-u_{2}$.
 Further, without losing generality
we assume that $C\geq0$ and for 
$$
w:=u-C(T  -t)
$$
 find that
$$
\delta _{\tau}w+L^{\alpha}_{h}w =
\delta _{\tau}u+L^{\alpha}_{h}
u +C+c^{\alpha}C(T-t)\geq\delta _{\tau}u+ L^{\alpha}_{h}u +C,
$$
$$
\delta _{\tau}w+ \sup_{\alpha\in A}L^{\alpha}_{h}w\geq0,\quad
w+\varepsilon\delta _{\tau}w
+\varepsilon \sup_{\alpha\in A}L^{\alpha}_{h}w\geq w
\quad\text{on}\quad Q.
$$
where $\varepsilon>0$ is any number. 

Next, looking at the proof of Lemma \ref{lemma 9.7.1} we
see that we can 
choose $\varepsilon$ so that, for $\gamma=1$ and
 any $\alpha\in A$ in \eqref{9.31.2} we have
$p^{\alpha}_{k}\geq0$, $p^{\alpha}\geq0$. 
Then for any function $\psi\geq w$ we have
\begin{equation}
                                                 \label{9.9.02}
\psi+\varepsilon\delta _{\tau}\psi+
\varepsilon \sup_{\alpha\in A}L^{\alpha}_{h}
\psi\geq w\quad\text{on}\quad Q.
\end{equation}

Take a rather small constant $\gamma>0$ 
to be specified later and take the function $\xi(t)$
from \eqref{9.31.1}.
Also introduce
$$
 \eta(x)=\cosh(\mu|x|),
\quad\zeta=\xi\eta,\quad
N_{0}=\sup_{\bar\cM_{T}}\frac{w_{+}}{\zeta}.
$$
Notice that by \eqref{9.9.2}
and by straightforward computations
$$
\sup_{\alpha\in A}L^{\alpha}_{h}\eta(x)
\leq\sup_{\alpha\in A}L^{\alpha}\eta(x)
$$
$$
+N_{1}(h^{2}+h)\cosh(\mu|x|+\mu K )\leq N_{2} \cosh(\mu|x|+\mu K ),
$$
where $N_{i}$ depend  only on $K$, $\mu$, and $d_{1}$.
It is seen as well that one can take $N_{2}$,
so that $N_{2}(K,d_{1},\mu)\to0$ as $\mu\downarrow0$
and $N_{2}(K,d_{1},0)=0$ even if $h>1$.
Also note that (cf.~\eqref{10.2.1})
$$
\delta _{\tau}\xi(t)=\xi(t) \tau ^{-1}
( \gamma -1).
$$
Therefore,
$$
\delta _{\tau}\zeta +
 \sup_{\alpha\in A}L^{\alpha}_{h}\zeta 
\leq\zeta (\tau^{-1}( \gamma -1)
+N_{3} )= \kappa\zeta,
$$
where
$$
N_{3}=N_{2} \sup_{x}\frac{\cosh(\mu|x|+\mu K)}
{\cosh(\mu|x|)}<\infty,\quad
\kappa=\kappa(\gamma):=\tau^{-1}( \gamma -1)
+N_{3}.
$$

Now set $\tau^{*}=N_{3}^{-1}$ and assume that
$\tau<\tau^{*}$. Upon noticing that $\kappa(0)<0$
and $\kappa(1)\geq0$ we see that  we can take
 $\gamma$ so   that $\kappa<0$ and $1+\kappa\varepsilon>0$.

After that for $\psi=N_{0}\zeta$ equation \eqref{9.9.02}
implies that
$$
N_{0}\zeta(1+\kappa\varepsilon)=
N_{0}\zeta+\kappa\varepsilon N_{0}\zeta\geq w
$$
on $Q$. Since the right-hand side is nonpositive
on $\bar\cM_{T}\setminus Q$, the inequality holds
on $\bar\cM_{T}$ and by the definition of $N_{0}$
implies that $N_{0}(1+\kappa\varepsilon)\geq N_{0}$.
By recalling that $\kappa<0$ we obtain
$N_{0}=0$, $w\leq0$ and \eqref{9.8.3} follows.

To prove the second assertion of the lemma it suffices
to add that if $\mu=0$, then $N_{3}=N_{2}=0$.
The lemma is proved.

Three completely standard
applications of the comparison principle follow.
\begin{corollary}
                                         \label{corollary 9.9.2}
Let a constant $c_{0}\geq0$ be such that
$$
\tau^{-1}(e^{c_{0}\tau}-1)\leq\lambda.
$$
Then
$$
|v_{\tau,h}(t,x)|\leq K\frac{1-e^{-\lambda
(T+\tau)}}
{\lambda}
+e^{-c_{0}(T -t)} \sup_{x}|g|
$$
 on $\bar H_{T}$ with natural interpretation
of this estimate if $\lambda=0$, that is
$$
|v_{\tau,h}|\leq K(T+\tau)+\sup_{x}|g|.
$$
 
\end{corollary}

To prove the corollary we observe that
it suffices to concentrate on $\bar\cM_{T}$.
Then we pass from $T$ to $T'$ thus reducing the general
case  to the one with $T=n\tau$, where $n$ is an integer. 
Next, define
$$
 N_{1}=\sup|g|,\quad
\xi(t)=K\lambda^{-1}(1-e^{-\lambda(T-t)})
+e^{-c_{0}(T-t)}N_{1} 
$$
if $\lambda>0$ with natural modification for $\lambda=0$.
We have $\xi\geq g=v_{\tau,h}$ on $\bar\cM_{T}
\setminus\cM_{T}$  whereas
 on $\cM_{T}$  
$$
\delta_{\tau}\xi(t)
-\lambda\xi(t)
=-K\bigg[e^{\lambda t-\lambda T}
\bigg(\frac   {e^{ \lambda\tau}-1}{\tau\lambda}-1\bigg)+1
\bigg]
$$
$$
+N_{1}\tau^{-1}(e^{c_{0}\tau}-1)
e^{-c_{0}(T-t)}
-\lambda N_{1}e^{-c_{0}(T-t)}\leq -K,
$$
so that
$$
\delta _{\tau}\xi + 
\sup_{\alpha\in A}[L^{\alpha}_{h}\xi +f ^{\alpha}]\leq0.
$$
By the lemma $v_{\tau,h}\leq\xi$ on $\bar\cM_{T}$.
Similarly one proves that $v_{\tau,h}\geq-\xi$.

\begin{corollary}
Let $u_{1}$ and $u_{2}$ be bounded solutions
of \eqref{8.19.1} in $H_{T}$ with terminal condition
$u_{1}(T,x)=g_{1}(x)$ and $u_{2}(T,x)=g_{2}(x)$,
where $g_{1}$ and $g_{2}$ are given bounded functions.
Then under the conditions of
Corollary \ref{corollary 9.9.2} we have
\begin{equation}
                                             \label{9.31.3}
u_{1}(t,x) \leq u_{2}(t,x)+e^{-c_{0}(T-t)}\sup(g_{1}-g_{2})_{+}
\end{equation}
in $\bar{H}_{T}$.
\end{corollary}

To prove this it suffices to replace $u_{2}$
in Lemma \ref{lemma 9.9.2} with 
the right-hand side of \eqref{9.31.3}.

\begin{corollary}
                                         \label{corollary 9.9.1}
Assume that there is a constant $R$ such that
 $f^{\alpha}(t,x)=g(x)=0$
if $|x|\geq R$. Then 
$$
\lim_{|x|\to\infty}\sup_{[0,T]}|v_{\tau,h}(t,x)|=0.
$$
\end{corollary}

For the proof take a unit $l\in\bR^{d}$ and
for small $\gamma\in(0,1)$ consider
$$
\zeta=\xi\eta, 
\quad\eta =e^{\gamma (x,l)},
$$
where $\xi$ is taken from the proof of the lemma.
It is a matter of very simple computations that
 $
L^{\alpha}_{h}\eta\leq N\gamma \eta,
 $
where $N$ is independent of $l$, $\gamma$, $\alpha$, and $t,x$.
It follows that
$$
\delta^{T}_{\tau}\zeta+ 
\sup_{\alpha\in A} L^{\alpha}_{h}\zeta\leq[\tau^{-1}(\gamma-1)+
N\gamma ]\zeta\leq0
$$
if $\gamma$ is sufficiently small. If needed we reduce further
the value of $\gamma$ to have  $\tau<\tau^*(K,d_{1},\gamma )$.
 Then on
$$
Q=\{(t,x)\in\cM_{T}:(x,l)\leq -R\} ,
$$
where $f^{\alpha}=0$,
we have
$$
\delta^{T}_{\tau}N\zeta+ 
\sup_{\alpha\in A}[L^{\alpha}_{h}N\zeta+f ^{\alpha}]\leq0
$$
for any constant $N>0$.  On
$\bar\cM_{T}\setminus Q$ it holds that 
$$
\zeta(t,x)\geq\exp(-\gamma R  )
\quad\text{if}\quad
t\in[0,T),
$$
$$
\zeta(t,x)\geq \exp(-\gamma |x|  )\quad\text{if}
\quad t= T,
$$
 which shows that
$N\zeta\geq v_{\tau,h}$ on
$\bar\cM_{T}\setminus Q$ for sufficiently large $N$.
By Lemma \ref{lemma 9.9.2} we obtain
$v_{\tau,h}\leq N\zeta$ in $\bar\cM_{T}$ 
and due to the arbitrariness
of $l$ we conclude  $$v_{\tau,h}\leq N\xi(0)
\exp(-\gamma |x| ).$$
Similarly, one proves that
$$v_{\tau,h}\geq -N\xi(0)
\exp(-\gamma |x| )
$$
and the result follows if we restrict ourselves
to considering $v_{\tau,h}$ only on $\bar\cM_{T}$.
But since every mesh \eqref{9.9.5} can be treated in the same way
and our constants stay the same,
we get the result as stated.

\begin{corollary}
                                         \label{corollary 9.9.3}
Let $h,\tau\leq K$.
Fix $(s_{0},x_{0})\in\bar\cM_{T}$ and set
$$
\nu=\sup_{(s_{0},x) \in\bar\cM_{T}}
\frac{|v_{\tau,h}(s_{0},x)-v_{\tau,h}(s_{0},x_{0}) 
|}{|x-x_{0}|}.
$$
Then for all $(t_{0},x_{0})\in\bar\cM_{T}$ with 
$s_{0}-1\leq t_{0}\leq s_{0}$ we have
$$
|v_{\tau,h}(s_{0},x_{0})-v_{\tau,h}(t_{0},x_{0})|
\leq N(\nu+1)|s_{0}-t_{0}|^{1/2},
$$
where $N$ depends only on $K$ and $d_{1}$.
\end{corollary}

To prove this we may assume that $s_{0}>0$.
Also, shifting the origin of the time axis
allows us to assume that $t_{0}=0$,
so that $s_{0}\leq1$. Then
 fix a constant $\gamma>0$,
define $s_{0}'$ as the least $n\tau$, $n=1,2,...$,
such that $s_{0}\leq n\tau$ and on $\bar\cM_{s_{0}}$ set
$$
\xi(t)=e^{s'_{0}-t}\quad t<s_{0} ,\quad
\xi(t)=1\quad t\geq s_{0} ,
$$
$$
\eta=|x-x_{0}|^{2},\quad
 \zeta=\xi\eta,
$$
$$
\psi =\gamma\nu
[\zeta+
\kappa(s_{0}-t)]+
K(s_{0}-t)+\gamma^{-1}\nu+v_{\tau,h}(s_{0},x_{0}),
$$
where $\kappa>0$ is a constant to be specified later.
It is easy to check that 
$\delta^{s_{0}}_{\tau}\xi=-\theta\xi$ on $Q=\cM_{s_{0}}$,
where
$$
\theta:=\tau^{-1}(1-e^{-\tau})\geq K^{-1}(1-e^{-K}).
$$
Also in $\cM_{s_{0}}$
$$
L^{\alpha}_{h}\eta(t,x)=2 a^{\alpha}_{k}(t,x)
|\ell_{k}|^{2}+b^{\alpha}_{k}(t,x)(\ell_{k},
2(x-x_{0})+h\ell_{k})
$$
$$
-c^{\alpha}(t,x)\eta(t,x)\leq N_{1}(1+|x-x_{0}|),
$$
$$
\delta^{s_{0}}_{\tau}\zeta(t,x)+
L^{\alpha}_{h}\zeta(t,x)\leq
[N_{1}(1+|x-x_{0}|)-\theta|x-x_{0}|^{2}]\xi(t)
$$
$$
\leq N_{2}(1+|x-x_{0}|)-\theta|x-x_{0}|^{2},
$$
where the constants $N_{i}$ depend only on
$K$  and $d_{1}$. It follows that
in~$\cM_{s_{0}}$
$$
\delta^{s_{0}}_{\tau}\psi +L^{\alpha}_{h}\psi 
+f^{\alpha} \leq\gamma\nu[
N_{2}(1+|x-x_{0}|)-\theta|x-x_{0}|^{2}-\kappa].
$$
As is easy to see there is $\kappa>0$
depending only on $N_{2}$ such that the right-hand side
is negative for all $x$.

Furthermore, 
$$
\psi(s_{0},x)=\nu(\gamma|x-x_{0}|^{2}+\gamma^{-1})
+v_{\tau,h}(s_{0},x_{0})
$$
$$
\geq\nu|x-x_{0}|+v_{\tau,h}(s_{0},x_{0})\geq v_{\tau,h}(s_{0},x).
$$
By Lemma \ref{lemma 9.9.2} applied to
$\cM_{s_{0}} $
in place of $\cM_{T}$ we conclude
$$
v_{\tau,h}(t,x_{0})\leq\psi(t,x_{0})
=\gamma\nu\kappa(s_{0}-t)+\gamma^{-1}\nu+K(s_{0}-t)
+v_{\tau,h}(s_{0},x_{0}).
$$
Minimizing with respect to $\gamma>0$ yields
$$
v_{\tau,h}(t,x_{0})-v_{\tau,h}(s_{0},x_{0})
\leq2\nu\kappa^{1/2}|s_{0}-t|^{1/2}+K
 s_{0} ^{1/2}|s_{0}-t|^{1/2}.
$$
Thus we  obtain a one-sided estimate of
$v_{\tau,h}(t,x_{0})-v_{\tau,h}(s_{0},x_{0})$.
The estimate from the other side is obtained similarly
by considering 
$$-
\gamma\nu
[\zeta+
\kappa(s_{0}-t)]-
K(s_{0}-t)-\gamma^{-1}\nu+v_{\tau,h}(s_{0},x_{0})
$$
in place of $\psi$.

One more simple consequence of Lemma \ref{lemma 9.7.1}
and Corollary \ref{corollary 9.9.2} is the following
stability result.
\begin{lemma}
                                              \label{lemma 9.9.3}
Let functions $f^{\alpha}_{n}$ and $g_{n}$, $n=1,2,...$,
satisfy the same conditions as $f^{\alpha},g$
with the same constants and let $v^{n}_{\tau,h}$
be the unique solutions of problems
\eqref{8.19.1}-\eqref{9.11.4} with 
$f^{\alpha}_{n}$ and $g_{n}$ in place 
of $f^{\alpha}$ and $g$, respectively.
Assume that on $\bar{H}_{T}$
$$
\lim_{n\to\infty}\sup_{\alpha\in A}(
|f^{\alpha}-f^{\alpha}_{n}|+|g-g_{n}|)=0.
$$
Then $v^{n}_{\tau,h}\to v_{\tau,h}$
on  $\bar{H}_{T}$.
\end{lemma}

Proof. It suffices again to concentrate on
$\bar\cM_{T}$ and observe that any subsequence
of uniformly bounded functions $v^{n}_{\tau,h}$
which converges at any point of $\bar\cM_{T}$
will converge to a solution of the original
problem \eqref{8.19.1}-\eqref{9.11.4},
which is unique and equals $v_{\tau,h}$.
Therefore, the whole sequence
converges to $v_{\tau,h}$.
The lemma is proved.

\mysection{Some technical tools}
                                       \label{section 10.2.3}

Set
$$
T_{h,l}u(x):=u(x+hl).
$$
\begin{lemma}
                                             \label{lemma 6.18.1}
For any functions $u(x),v(x)$, $h>0$,
and $l\in\bR^{d}$ we have
\begin{equation}
                                                \label{10.2.4}
T_{h,-l}T_{h,l}u=u,
\end{equation}
\begin{equation}
                                                \label{6.21.1}
T_{h,l}\delta_{h,-l}=\delta_{h,-l}T_{h,l}=-\delta_{h,l},\quad
T_{h,-l}\delta_{h,l}=\delta_{h,l}T_{h,-l}=-\delta_{h,-l},
\end{equation}
\begin{equation}
                                                 \label{6.18.1}
\delta_{h,l}(uv)=(\delta_{h,l}u)v+(T_{h,l}u)\delta_{h,l}v=
v\delta_{h,l}u+u\delta_{h,l}v
+h(\delta_{h,l}u)\delta_{h,l}v,
\end{equation}
\begin{equation}
                                               \label{6.18.2}
\Delta_{h,l}(uv)=v\Delta_{h,l}u+u\Delta_{h,l}v
+(\delta_{h,l}u)\delta_{h,l}v+(\delta_{h,-l}u)\delta_{h,-l}v.
\end{equation}
In particular,
\begin{equation}
                                               \label{6.18.3}
\Delta_{h,l}(u^{2})=2u\Delta_{h,l}u
+(\delta_{h,l}u)^{2}+(\delta_{h,-l}u)^{2}.
\end{equation}
\end{lemma}

Proof. Equations \eqref{10.2.4}, \eqref{6.21.1},
and  \eqref{6.18.1}  are almost trivial. They yield
equation (\ref{6.18.3}) because
$$
-\Delta_{h,l}(u^{2})=\delta_{h,-l}[(\delta_{h,l}u)u]
+\delta_{h,-l}[(T_{h,l}u) \delta_{h,l}u]
$$
$$
=[-(\Delta_{h,l}u)u+(T_{h,-l}\delta_{h,l}u)\delta_{h,-l}u]
$$
$$
+[(\delta_{h,-l}T_{h,l}u)\delta_{h,l}u-
(T_{h,-l}T_{h,l}u)\Delta_{h,l}u].
$$

Equation (\ref{6.18.2}) is obtained by
polarizing (\ref{6.18.3}) that is by comparing
the coefficient of $\lambda$ in (\ref{6.18.3})
applied to $u+\lambda v$ in place of $u$.
The lemma is proved.

\begin{lemma}
                                               \label{lemma 6.18.2}

Let $u,v,w$ be functions on $\bR^{d}$, $l,x_{0}\in\bR^{d}$,
$h>0$. Assume that $v(x_{0})\leq0$ and $w(x_{0})\leq0$.
Then at $x_{0}$ it holds that
\begin{equation}
                                             \label{6.18.5}
-\delta_{h,l}v\leq\delta_{h,l}(v_{-}),\quad
-\Delta_{h,l}v\leq\Delta_{h,l}(v_{-}),
\end{equation}
\begin{equation}
                                                   \label{9.23.5}
-\delta_{h,l}(u_{-})\leq
[\delta_{h,l}((u+v)_{-})]_{-}+[\delta_{h,l}(v_{-})]_{+},
\end{equation}
$$
(\Delta_{h,l}u)_{-}\leq 
[\delta_{h,-l}((\delta_{h,l}u+v)_{-})]_{-}+
 [\delta_{h,l}((\delta_{h,-l}u+w)_{-})]_{-}
$$
\begin{equation}
                                                   \label{9.23.6}
+[\delta_{h,-l}(v_{-})]_{+}+[\delta_{h, l}(w_{-})]_{+},
\end{equation}
\begin{equation}
                                          \label{6.18.6}
|\Delta_{h,l}u|\leq 
|\delta_{h,-l}((\delta_{h,l}u )_{-})|+
 |\delta_{h,l}((\delta_{h,-l}u )_{-})| ,
\end{equation}
\begin{equation}
                                             \label{6.21.2}
 |\Delta_{h,l}u|\leq 
|\delta_{h,-l}((\delta_{h,l}u )_{+})|+
 |\delta_{h,l}((\delta_{h,-l}u )_{+})|.
\end{equation}

\end{lemma}

Proof. We use the formulas $-\alpha\leq \alpha_{-}$ 
and $v(x_{0})=-v_{-}(x_{0})$
and get
$$
-h\delta_{h,l}v(x_{0})
=v(x_{0})-v(x_{0}+hl) 
\leq
-v_{-}(x_{0})+v_{-}(x_{0}+hl),
$$
which is the first inequality in (\ref{6.18.5}). 
The second one is obtained by summing up the first inequality
corresponding to $l$ and $-l$.
 
 While proving \eqref{9.23.5} we may assume
that $u(x_{0})<0$ since otherwise the
left-hand side   is negative. In that case by noting that
by subadditivity:$(\alpha+\beta)_{-}\leq
\alpha_{-}+\beta_{-}$, we have
$$
 -u_{-}
\leq  -( u+v)_{-}+
 v_{-},\quad
 -T_{h,l}u_{-}
\leq  -T_{h,l}( u+v)_{-}+
 T_{h,l}v_{-}
$$
everywhere, whereas since $u(x_{0})\leq0,v(x_{0})\leq0$,
we have at $x_{0}$
$$
u_{-}=(u+v)_{-}-v_{-}.
$$
We conclude  that at $x_{0}$
$$
-\delta_{h,l}u_{-}\leq-\delta_{h,l}(u+v)_{-}+\delta_{h,l}v_{-}
$$
and   \eqref{9.23.5} follows.

In the proof of \eqref{9.23.6} we may assume that
$\Delta_{h,l}u(x_{0})\leq0$. Then 
owing to (\ref{6.21.1})  at $x_{0}$
$$
(\Delta_{h,l}u)_{-}= \delta_{h,-l} \delta_{h,l}u 
= \delta_{h,-l}((\delta_{h,l}u)_{+})
 -\delta_{h,-l}((\delta_{h,l}u)_{-})
$$
$$
=T_{h,l}\delta_{h,-l}((-\delta_{h,-l}u)_{+})
-\delta_{h,-l}((\delta_{h,l}u)_{-})
$$
$$
=-\delta_{h,l}((\delta_{h,-l}u)_{-})
-\delta_{h,-l}((\delta_{h,l}u)_{-}).
$$
This and \eqref{9.23.5} imply (\ref{6.18.6}) .

If $\Delta_{h,l}u(x_{0})\leq0$, \eqref{6.18.6}
follows from \eqref{9.23.6} with $v\equiv w\equiv0$.
Therefore, we may concentrate on the case that 
 $\Delta_{h,l}u(x_{0})\geq0$.
By applying \eqref{9.23.6} with $v\equiv w\equiv0$ to $-u$ in place of $u$
and using (\ref{6.21.1}) we get at $x_{0}$ that
$$
|\Delta_{h,l}u|\leq
|\delta_{h,-l}((-\delta_{h,l}u)_{-})|+
 |\delta_{h,l}((-\delta_{h,-l}u)_{-})|
$$
$$
=|\delta_{h,-l}(T_{h,l}(\delta_{h,-l}u)_{-})|+
 |\delta_{h,l}(T_{h,-l}(\delta_{h,l}u)_{-})|
$$
$$
=|T_{h,l}\delta_{h,-l}((\delta_{h,-l}u)_{-})|+
 |T_{h,-l}\delta_{h,l}((\delta_{h,l}u)_{-})|
$$
$$
=|\delta_{h, l}(( \delta_{h,-l}u)_{-})|+
 |\delta_{h,-l}(( \delta_{h, l}u)_{-})|.
$$
This proves (\ref{6.18.6}).

Equation (\ref{6.21.2}) is obtained from (\ref{6.18.6})
by substituting $-u$ in place of $u$.
The lemma is proved.

\mysection{Main estimates}
                                       \label{section 10.2.4}

We take $\tau,h,T$, and $\cM_{T}$ from Section \ref{section 9.10.1},
fix an $\varepsilon\in[0,Kh]$ and a unit vector 
$l\in\bR^{d}$ and introduce
$$
\bar\cM_{T}(\varepsilon):=\{(t,x+i\varepsilon l):(t,x)\in
\bar\cM_{T} ,
i=0,\pm1,...\}.
$$
 Let $Q\subset\bar\cM_{T}(\varepsilon)$
 be a nonempty {\em finite\/} set
 and $u$   a   function on $\bar\cM_{T}(\varepsilon)$
satisfying \eqref{8.19.1} in $Q':=Q\cap([0,T)\times\bR^{d})$.

Set
$$
 Q^{o}_{\varepsilon} =\{(t,x)\in Q': 
(t+\tau_{T}(t),x), (t,x\pm h\ell_{k}),
(t,x\pm\varepsilon l) 
\in Q,\forall k= 1,...,  d_{1} \},
$$
$$
 \partial_{\varepsilon} Q=Q\setminus Q^{o}_{\varepsilon} .
$$

Instead of Assumption \ref{assumption 8.25.1}
in this section we use  the following.

\begin{assumption}
                                   \label{assumption 8.25.2}

 For 
$$
\psi 
=  b^{\alpha}_{k},c^{\alpha}
-\lambda,
f^{\alpha}  ,\quad k=\pm1,...,\pm d_{1},\quad
 \alpha\in A
$$
 we have in $Q^{o}_{\varepsilon}$ that
\begin{equation}
                                                 \label{8.24.03}
|\psi  |\leq K,\quad|\delta_{h,\ell_{k}}\psi |,
|\delta_{\varepsilon,\pm l}\psi |\leq K,
\quad b^{\alpha}_{k},c^{\alpha}
-\lambda,\lambda\geq0, 
\end{equation}
\begin{equation}
                                                 \label{8.25.3}
0\leq a^{\alpha}_{k}  \leq K,\quad
|\delta_{h,\ell_{k}}a^{\alpha}_{k} |,
|\delta_{\varepsilon,\pm l}a^{\alpha}_{k} |\leq
K \sqrt{a^{\alpha}_{k} }+Kh .
\end{equation}
\end{assumption}

\begin{theorem}
                                            \label{theorem 8.19.1}
There is a constant $N\in(0,\infty)$
 depending only on $K$ and
  $d_{1}$,   such that if for a number
$c_{0}\geq0$ it holds that
\begin{equation}
                                                \label{9.11.11}
\lambda+\frac{1-e^{-c_{0}\tau}}{\tau}> N,
\end{equation}
then for $\varepsilon\in(0,Kh]$ on $Q$
\begin{equation}
                                                  \label{8.19.4}
|\delta_{\varepsilon,\pm l}u| \leq 
Ne^{c_{0}(T+\tau)}\big(1+|u|_{0,Q}
+\max_{\partial_{\varepsilon} Q}( \max_{k}
|\delta_{h ,\ell_{k}}u|+|\delta_{\varepsilon, l}u|
+|\delta_{\varepsilon,-l}u|) \big).
\end{equation}
\end{theorem}

 Before proving the theorem
we do some preparations. Denote
$$
 h_{k}=h,\quad
k=\pm1,...,\pm d_{1},\quad
h_{\pm(d_{1}+1)}=\varepsilon,\quad\ell_{\pm(d_{1}+1)}=\pm l, 
$$
and let $ r $ be an index   running 
through $\{ \pm1,...,\pm(d_{1}+1)\}$
and $k $ through $\{ \pm1,...,\pm d_{1}\}$.
 
Take a constant $c_{0}\geq0$ and introduce $T'$
as the least $n\tau$, $n=1,2,...$, such that $n\tau\geq T$,
$$
\xi(t)=e^{c_{0}t},\quad t<T,\quad\xi(T)=e^{c_{0}T'},
$$
$$
v =\xi u,\quad
v_{r }  =\delta_{h_{r },\ell_{r }}v,\quad
v_{r }^{-}=(v_{r })_{-},
$$
$$
M_{0}=\max_{Q}|v|,\quad
M_{1}=\max_{Q,r }|v_{r }|.
$$

Let $(t_{0},x_{0})$ be a point in $Q$ at which
$$
V:=\sum_{r} (v_{r }^{-})^{2}
$$
attains its maximum value in $Q$. 

Observe that for each $(t,x)\in Q^{o}_{\varepsilon} $ and $r $ we have
$$
(t,x+ h_{r }\ell_{r })\in Q  
$$
 and 
$$
\text{either}\quad v_{r }
(t,x)\leq0\quad\text{or}\quad -v_{r }(t,x)=
 v_{-r }(t,x+h_{r }\ell_{r })\leq0.
$$
In the first case 
$$
|v_{r }(t,x)|\leq V^{1/2}(t,x)\leq V^{1/2}(t_{0},x_{0}),
$$
whereas in the second case
$$
|v_{r }(t,x)|\leq V^{1/2}(t,x+h_{r }\ell_{r })
\leq V^{1/2}(t_{0},x_{0}).
$$
It follows that
\begin{equation}
                                               \label{8.22.1}
M_{1}\leq  \max_{\partial_{\varepsilon} Q,r }|v_{r }|
+V^{1/2}(t_{0},x_{0}),
\end{equation}
\begin{equation}
                                               \label{8.19.5}
|\delta_{\varepsilon,\pm l}u|
\leq e^{c_{0}T'} \max_{\partial_{\varepsilon} Q,r }|
\delta_{h_{r },\ell_{r }}u|+V^{1/2}(t_{0},x_{0})
\end{equation}
on $Q$
and we need only estimate $V^{1/2}(t_{0},x_{0})$.

Furthermore,  obviously
$$
V^{1/2}(t ,x )\leq2d_{1}
\max_{r }|v_{r }(t ,x )|
\leq2d_{1}e^{c_{0}(T+\tau)}
\max_{r }|\delta_{h_{r },\ell_{r }}u(t ,x )|,
$$
so that while estimating $V^{1/2}(t_{0},x_{0})$
we may assume that
\begin{equation}
                                               \label{10.3.2}
(t_{0},x_{0})\in Q^{o}_{\varepsilon}.
\end{equation}

Notice that there is a sequence $\alpha_{n}\in A$
such that
$$
\delta^{T}_{\tau}u(t_{0},x_{0})+
\lim_{n\to\infty}[a^{\alpha_{n}}_{k}(t_{0},x_{0})
\Delta_{h,\ell_{k}}u(t_{0},x_{0} )
+b^{\alpha_{n}}_{k}(t_{0},x_{0})
\delta_{h,\ell_{k}}u(t_{0},x_{0})
$$
$$
-c^{\alpha_{n}}(t_{0},x_{0})u(t_{0},x_{0})+
f^{\alpha_{n}}(t_{0},x_{0})]
=\delta^{T}_{\tau}u(t_{0},x_{0})
$$
$$
+F(\Delta_{h,\ell_{k}}u(t_{0},x_{0} ),
\delta_{h,\ell_{k}}u(t_{0},x_{0}),u(t_{0},x_{0}),t_{0},x_{0} )=0.
$$
Owing to Assumption \ref{assumption 8.25.2}
there is a subsequence $\{n'\}\subset\{1,2,...\}$
and  functions
$\bar{a}_{k}
(t,x),\bar{b}_{k}(t,x),\bar{c}(t,x),\bar{f}(t,x)$
such that they satisfy 
Assumption \ref{assumption 8.25.2} changed
in an obvious way and
$$
(\bar{a}_{k}
(t,x),\bar{b}_{k}(t,x),\bar{c}(t,x),\bar{f}(t,x))
$$
$$
=\lim_{n'\to\infty}
(a^{\alpha_{n'}}_{k}(t,x),b^{\alpha_{n'}}_{k}(t,x),
c^{\alpha_{n'}} (t,x),f^{\alpha_{n'}} (t,x))
$$
on $Q$. Obviously,  at $(t_{0},x_{0})$ we have
\begin{equation}
                                               \label{9.7.1}
\delta^{T}_{\tau}u+\bar{a}_{k}\Delta_{h_{k},\ell_{k}}u+\bar{b}_{k}
\delta_{h_{k},\ell_{k}}u-\bar{c}u+\bar{f}=0,
\end{equation}
and for any $r$ ($ =\pm1,...,\pm (d_{1}+1)$) 
owing to \eqref{10.3.2}
\begin{equation}
                                               \label{9.7.2}
T_{h_{r },\ell_{r }}[\delta^{T}_{\tau}u+
\bar{a}_{k}\Delta_{h_{k},\ell_{k}}u+\bar{b}_{k}
\delta_{h_{k},\ell_{k}}u-\bar{c}u+\bar{f}]\leq0,
\end{equation}
where and below 
for simplicity of notation
we drop $(t_{0},x_{0})$ in the arguments of functions
which we are dealing with.

\begin{lemma} 
                                            \label{lemma 8.21.1}
For all $k=\pm1,...,\pm d_{1}$ 
  at $(t_{0},x_{0})$ we have
\begin{equation}
                                               \label{6.22.1}
 v_{r }^{-}
\Delta_{h_{k},\ell_{k}}v_{r }  
\geq0.
\end{equation}
Furthermore, there is a constant $N\in(0,\infty)$
 depending only on $K$
and  $d_{1}$,   such that at $(t_{0},x_{0})$
$$
\hat{\lambda}V + (1/2)v_{r }^{-}
\bar{a}_{k}\Delta_{h_{k},\ell_{k}}v_{r }
+(1/2)I
+v_{r }^{-}
(\delta_{h_{r },\ell_{r }}\bar{a}_{k})
\Delta_{h_{k},\ell_{k}}v
$$
\begin{equation}
                                                \label{6.20.2}
+v_{r }^{-} 
h_{r }(\delta_{h_{r },\ell_{r }}
\bar{a}_{k})\Delta_{h_{k},\ell_{k}}v_{r }
\leq N(e^{c_{0} T'}+M_{0}+M_{1})M_{1},
\end{equation}
where
$$
\hat{\lambda}=\lambda+\frac{1-e^{-c_{0}\tau}}{\tau},\quad
I=
 \sum_{r }\bar{a}_{k}  (\delta_{h_{k},\ell_{k}}
v_{r }^{-})^{2}.
$$
\end{lemma}

Proof. By Lemma \ref{lemma 6.18.1} and Lemma \ref{lemma 6.18.2}
(with $v_{r }$ in place of $v$) 
$$
0\geq\Delta_{h_{k},\ell_{k}}\sum_{r }(v_{r }^{-})^{2}
=2v_{r }^{-}
\Delta_{h_{k},\ell_{k}}v_{r }^{-}
+\sum_{r }[(\delta_{h_{k},\ell_{k}}v_{r }^{-})^{2}
+(\delta_{h_{k},\ell_{-k}}v_{r }^{-})^{2}]
$$
$$
\geq-2v_{r }^{-}
\Delta_{h_{k},\ell_{k}}v_{r } 
+\sum_{r }[(\delta_{h_{k},\ell_{k}}v_{r }^{-})^{2}
+(\delta_{h_{k}\ell_{-k}}v_{r }^{-})^{2}].
$$
This obviously yields (\ref{6.22.1}) and also that
$$
I\leq v_{r }^{-}
\bar{a}_{k}\Delta_{h_{k},\ell_{k}}v_{r },
$$
which in turn implies that to prove (\ref{6.20.2})
it suffices to prove that
$$
\hat{\lambda}V +  v_{r }^{-}
\bar{a}_{k}\Delta_{h_{k},\ell_{k}}v_{r }
+v_{r }^{-}
(\delta_{h_{r },\ell_{r }}\bar{a}_{k})
\Delta_{h_{k},\ell_{k}}v
$$
\begin{equation}
                                                \label{6.21.01}
+v_{r }^{-} 
h_{r }(\delta_{h_{r },\ell_{r }}
\bar{a}_{k})\Delta_{h_{k},\ell_{k}}v_{r }
\leq N(e^{c_{0} T'}+M_{0}+M_{1})M_{1}.
\end{equation}

By subtracting  the inequalities 
\eqref{9.7.1} and \eqref{9.7.2} and using
(\ref{6.18.1}) we find 
\begin{equation}
                                                \label{6.18.7}
\delta^{T}_{\tau}(\xi^{-1}v_{r })+
\xi^{-1}\big[\bar{a}_{k}\Delta_{h_{k},\ell_{k}}v_{r }
+I_{1r }
+I_{2r }
+I_{3r }
+I_{4r }\big]\leq0,
\end{equation}
where (no summation in $r $)
$$
I_{1r }=(\delta_{h_{r },\ell_{r }}\bar{a}_{k})
\Delta_{h_{k},\ell_{k}}v,
$$
$$
I_{2r }=h_{r }(\delta_{h_{r },\ell_{r }}
\bar{a}_{k})\Delta_{h_{k},\ell_{k}}v_{r },
$$
$$
I_{3r }=(T_{h_{r },\ell_{r }}\bar{b}_{k})\delta_{h_{k},\ell_{k}}
v_{r }
+(\delta_{h_{r },\ell_{r }}\bar{b}_{k})\delta_{h_{k},\ell_{k}}v,
$$
$$
I_{4r }=
-(\delta_{h_{r },\ell_{r }}\bar{c})v-(T_{h_{r },\ell_{r }}\bar{c})
v_{r }+\xi\delta_{h_{r },\ell_{r }}\bar{f}.
$$

We multiply (\ref{6.18.7}) by $\xi v_{r }^{-}$
and sum up with respect to $r $. 

Observe that in $I_{4r }$
$$
\delta_{h_{r },\ell_{r }}\bar{f}
\geq-K,\quad|\delta_{h_{r },\ell_{r }}\bar{c}|
\leq K, 
$$
$$
- v_{r }^{-}(T_{h_{r },\ell_{r }}\bar{c})
 v_{r }= (T_{h_{r },\ell_{r }}\bar{c})
[( v_{r }^{-}]^{2}
\geq \lambda\sum_{r }[ v_{r }^{-}]^{2}
= \lambda V,
$$
since $\bar{c}\geq \lambda$. Therefore,
$$
 v_{r }^{-}I_{4r }
\geq-KM_{1}(e^{c_{0}T'}+M_{0})+\lambda V.
$$

By using the fact that $V$ attains its maximum in $Q$
at
$(t_{0},x_{0})\in Q^{0}_{\varepsilon}$ and using Lemma \ref{lemma 6.18.2}
(with $v_{r }$ in place of $v$) we get
$$
0\geq\delta_{h_{k},\ell_{k}}\sum_{r }(v_{r }^{-})^{2}
=2v_{r }^{-}\delta_{h_{k},\ell_{k}}
v_{r }^{-}
+\sum_{r }h_{k}
(\delta_{h_{k},\ell_{k}}
v_{r }^{-})^{2}
$$
$$
\geq2v_{r }^{-}\delta_{h_{k},\ell_{k}}
v_{r }^{-}
\geq-
2v_{r }^{-}\delta_{h_{k},\ell_{k}}
v_{r }.
$$
This result and the inequalities
$b_{k}\geq0$, $|\delta_{h_{r },\ell_{r }}\bar{b}_{k}|
\leq K$  yield
$$
-v_{r }^{-}(T_{h_{r },\ell_{r }}\bar{b}_{k})
\delta_{h_{k},\ell_{k}}
v_{r }\leq0,\quad
v_{r }^{-}I_{3r }
\geq-NM_{1}^{2}.
$$

Similarly,
$$
0\leq-\delta^{T}_{\tau} \sum_{r}(v_{r }^{-})^{2}
 \leq
2v_{r }^{-}\delta^{T}_{\tau}
v_{r },
$$
which implies that
$$
\xi v_{r }^{-}\delta^{T}_{\tau}(\xi^{-1}v_{r })
=\xi v_{r }^{-}[\xi^{-1}(t_{0}+\tau_{T}(t_{0}))\delta^{T}_{\tau}v_{r }
+v_{r }\delta^{T}_{\tau}\xi^{-1}]
$$
$$
=e^{-c_{0}\tau}v_{r }^{-}\delta^{T}_{\tau}v_{r }
-V\xi\delta^{T}_{\tau}\xi^{-1}\geq-V\xi\delta^{T}_{\tau}\xi^{-1}
= V\frac{1}{\tau}[1-e^{-c_{0}\tau}].
$$
 
By combining the above estimates   we come to (\ref{6.21.01})
and the lemma is proved.

{\bf Proof of Theorem \ref{theorem 8.19.1}}.
By Lemma \ref{lemma 8.21.1}
\begin{equation}
                                                  \label{6.19.1}
\hat{\lambda}V\leq N(e^{c_{0}(T+\tau)}+M_{0}+M_{1})M_{1}+J_{1}+J_{2},
\end{equation}
where
$$
J_{1}:= v_{r }^{-}
|(\delta_{h_{r },\ell_{r }}\bar{a}_{k})\Delta_{h_{k},\ell_{k}}v|
-(1/4)\sum_{r }
\bar{a}_{k}(\delta_{h_{k},\ell_{k}}
v_{r }^{-})^{2},
$$
$$
J_{2}:=J_{3}
-(1/2)\bar{a}_{k}
 v_{r }^{-}\Delta_{h_{k},\ell_{k}}v_{r } 
-(1/4)\sum_{r}
\bar{a}_{k}(\delta_{h_{k},\ell_{k}}v_{r }^{-})^{2},
$$
$$
J_{3}:=h_{r}v_{r }^{-}|(\delta_{h_{r },\ell_{r }}
\bar{a}_{k})\Delta_{h_{k},\ell_{k}}v_{r }|.
$$

First we estimate $J_{1}$. By Lemma \ref{lemma 6.18.2}
$$
|\Delta_{h_{k},\ell_{k}}v|\leq\sum_{r }|\delta_{h_{k},\ell_k}
v_{r }^{-}|+ 
\sum_{r }|\delta_{h_{k},\ell_{-k}}
v_{r }^{-}|.
$$
Also we recall (\ref{8.25.3}) and use the inequality
$$
h |\Delta_{h_{k},\ell_{k}}v|\leq 2M_{1}.
$$
Then we obtain
 $$
 v_{r }^{-}
|(\delta_{h_{r },\ell_{r }}\bar{a}_{k})\Delta_{h_{k},\ell_{k}}v|
\leq NM_{1}|
(\sqrt{\bar{a}_{k}}+h)\Delta_{h_{k},\ell_{k}}v|
$$
$$
\leq N M^{2}_{1}+ (1/4)
\sum_{r }\bar{a}_{k}(\delta_{h_{k},\ell_{k}}v_{r }^{-})^{2},
\quad
 J_{1}\leq N M^{2}_{1}.
$$

To estimate $J_{3}$ observe that 
$$
h_{r}\leq Kh,\quad |a|=2a_{-}+a,\quad h^{2}|
\Delta_{h_{k},\ell_{k}}v_{r }|
\leq4M_{1},
$$ 
so that
$$
J_{3}
\leq
N_{1}v_{r }^{-} h\sqrt{\bar{a}_{k}}
|\Delta_{h_{k},\ell_{k}}v_{r }|
+K^{2}v_{r }^{-}h^{2}| \sum_{k}
\Delta_{h_{k},\ell_{k}}v_{r }|
$$
$$
\leq N_{1}
v_{r }^{-} h\sqrt{\bar{a}_{k}} 
|\Delta_{h_{k},\ell_{k}}v_{r }|+
N_{2} M^{2}_{1}
=2N_{1}v_{r }^{-}h \sqrt{\bar{a}_{k}} 
(\Delta_{h_{k},\ell_{k}}v_{r })_{-}
$$
$$
+N_{1}v_{r }^{-} h\sqrt{\bar{a}_{k}} 
\Delta_{h_{k},\ell_{k}}v_{r }+N_{2} M^{2}_{1}.
$$
Here the summation in $r $ can be restricted
to $r $ such that
$$
v_{r }^{-}\ne0,
$$
when
by Lemma \ref{lemma 6.18.2} it holds that
$$
h(\Delta_{h_{k},\ell_{k}}v_{r })_{-}
\leq h|\Delta_{h_{k},\ell_{k}}(v_{r }^{-})|=
|(\delta_{h_{k},\ell_{k}}+\delta_{h_{k}\ell_{-k}})
(v_{r }^{-})|
$$
$$
\leq| \delta_{h_{k},\ell_{k}}(v_{r }^{-})|
+|\delta_{h_{k}\ell_{-k}}
(v_{r }^{-})|.
$$
Therefore,
$$
J_{3}\leq N_{1}v_{r }^{-} h\sqrt{\bar{a}_{k}} 
\Delta_{h_{k},\ell_{k}}v_{r }
+N_{2} M^{2}_{1}
+N_{3}M_{1}\bigg[
\sum_{r }
\bar{a}_{k}(\delta_{h_{k},\ell_{k}}
v_{r }^{-})^{2}
\bigg]^{1/2},
$$
$$
J_{2}\leq NM_{1}^{2} 
-(1/2)(\bar{a}_{k}-2N_{1}h\sqrt{\bar{a}_{k}})
 v_{r }^{-}\Delta_{h_{k},\ell_{k}}v_{r } .
$$

Finally, let 
$$
\cK=\{k:\bar{a}_{k}-2N_{1}h\sqrt{\bar{a}_{k}}\geq0\}.
$$
Then, for $k\not\in\cK$ we have 
$$
\sqrt{\bar{a}_{k}}\leq 2N_{1}h \quad
 \bar{a}_{k}\leq 4N_{1}^{2}h^{2},\quad
 |\bar{a}_{k}-2N_{1}h\sqrt{\bar{a}_{k}}|\leq Nh^{2}
$$
and by using (\ref{6.22.1}) and using again
the fact that
 $h^{2}|\Delta_{h_{k}l}\phi|\leq 4\sup|\phi|$ we conclude
that
$$
-(1/2)(\bar{a}_{k}-2N_{1}h\sqrt{\bar{a}_{k}})
 v_{r }^{-}
\Delta_{h_{k},\ell_{k}}v_{r }
$$
$$
\leq -(1/2)\sum_{k\in\cK}(\bar{a}_{k}-2N_{1}h\sqrt{\bar{a}_{k}})
 v_{r }^{-}\Delta_{h_{k},\ell_{k}}v_{r }+N  M^{2}_{1}
 \leq N  M^{2}_{1},
$$
$$
J_{2}\leq N  M^{2}_{1}.
$$
Coming back to (\ref{6.19.1}) we get 
$$
\hat{\lambda}V\leq N(e^{c_{0}(T+\tau)}+M_{0} +M_{1})M_{1},
$$
which due to (\ref{8.22.1}) leads to
\begin{equation}
                                                  \label{8.22.2}
\hat{\lambda}V\leq N(e^{c_{0}(T+\tau)}+M_{0} +\mu+V^{1/2})
(\mu+V^{1/2}),
\end{equation}
where 
$$
\mu:=\sup_{\partial_{\varepsilon}Q,r }|v_{r }|
\leq e^{c_{0}(T+\tau)}\sup_{\partial_{\varepsilon}Q,r }|
\delta_{h_{r },\ell_{r }}u|=:
e^{c_{0}(T+\tau)}\bar{\mu}.
$$
Also introduce
$$
\bar{M}_{0}=|u|_{0,Q},
\quad
\bar{V}=e^{-2c_{0}(T+\tau)}V
$$
and notice that 
$$
M_{0}\leq e^{c_{0}(T+\tau)}\bar{M}_{0}.
$$
Then (\ref{8.22.2}) yields
$$
\hat{\lambda}\bar{V}
\leq N(1+\bar{M}_{0} +\bar{\mu}+\bar{V}^{1/2})
(\bar{\mu}+\bar{V}^{1/2})
$$
$$
\leq N^{*}(1+\bar{M}_{0}^{2}+\bar{\mu}^{2}
+\bar{V}).
$$
If $\hat{\lambda}\geq N^{*}+1$, then
 we conclude that
$$
\bar{V}\leq N^{*}(1+\bar{M}_{0}^{2}+\bar{\mu}^{2}),
$$
which along with (\ref{8.19.5}) yields (\ref{8.19.4})
and proves the theorem.

The following theorem bears on estimates
of how close two solutions
of the Bellman finite-difference equations are
if the coefficients are close. It is a generalization
of Theorem \ref{theorem 8.19.1}.

In the rest of the section 
we take some objects $\hat{\sigma}^{\alpha}_{k},
\hat{b}^{\alpha}_{k},\hat{c}^{\alpha},\hat{\lambda},
\hat{f}^{\alpha}$ defined on $A\times[0,T]\times\bR^{d}$
and having the same sense as in Section \ref{section 10.4.1}.
We set 
$$
\hat{a}^{\alpha}_{k}=(1/2)|\hat{\sigma}^{\alpha}_{k}|^{2}.
$$ 
\begin{assumption}
                                    \label{assumption 8.25.3}
We have a finite set $Q\subset\bar\cM=\bar\cM(0)$ and
not only $a^{\alpha}_{k},
b^{\alpha}_{k},c^{\alpha},  \lambda,
f^{\alpha}$ satisfy Assumption \ref{assumption 8.25.2}
with $\varepsilon=0$
but $\hat{a}^{\alpha}_{k},
\hat{b}^{\alpha}_{k},\hat{c}^{\alpha},\hat{\lambda},
\hat{f}^{\alpha}$ 
satisfy Assumption \ref{assumption 8.25.2} 
with $\varepsilon=0$ as well. Moreover,
$\lambda=\hat{\lambda}$.

\end{assumption}

\begin{theorem}
                                            \label{theorem 8.24.1}
  Let
$u$  be a   function on $\bar\cM_{T}$
satisfying   \eqref{8.19.1}  in $Q\cap([0,T)\times\bR^{d})$ and
let   $\hat{u}$ be a function  on $\bar\cM_{T}$
  satisfying equation
 \eqref{8.19.1}  in $Q\cap([0,T)\times\bR^{d})$ with 
$\hat{a}^{\alpha}_{k},
\hat b^{\alpha}_{k},\hat c^{\alpha}, 
\hat f^{\alpha}$   in place of
$a^{\alpha}_{k},
b^{\alpha}_{k},c^{\alpha},  
f^{\alpha}$, respectively.
 
Assume that there is an $\varepsilon\in(0,Kh]$
such that for $k=\pm1,...,\pm d_{1}$ on $ Q^{o}_{0}$ we have
\begin{equation}
                                                \label{8.26.1}
|b^{\alpha}_{k}-\hat{b}^{\alpha}_{k}|
+|c^{\alpha} -\hat{c}^{\alpha} |
+|f^{\alpha} -\hat{f}^{\alpha} |\leq
K\varepsilon,
\end{equation}
\begin{equation}
                                                \label{8.26.2}
|a^{\alpha}_{k}-\hat{a}^{\alpha}_{k}|\leq
K\varepsilon
\sqrt{a^{\alpha}_{k}\wedge\hat{a}^{\alpha}_{k}}
+K\varepsilon h.
\end{equation}

We assert that there is a constant $N\in(0,\infty)$
 depending only on $K$ and 
  $d_{1}$,   such that if for a number
$c_{0}\geq0$ equation \eqref{9.11.11} holds,
then in $Q$
$$
| u-\hat{u}| \leq 
N\varepsilon e^{c_{0}(T+\tau)}\big(1+|u|_{0,Q}+|\hat u|_{0,Q}
$$
\begin{equation}
                                                  \label{8.25.1}
+\sup_{\partial_{0}Q}(\max_{k}|\delta_{h ,\ell_{k}}u|
+\max_{k}|\delta_{h ,\ell_{k}}\hat{u}|
+\varepsilon^{-1}|u-\hat{u}|)\big).
\end{equation}
\end{theorem}

Proof. 
We want to apply Theorem \ref{theorem 8.19.1}
to appropriate objects.
Consider $\bR^{d}$ as a subspace
of 
$$
\bR^{d+1}=\{x=(x',x^{d+1}):x'\in\bR^{d},
x^{d+1}\in\bR\}.
$$
 Take an integer $m\geq1/\varepsilon$ and
introduce $l=(0,...,0,1)\in\bR^{d+1}$. Then
$$
\bar\cM_{T}(\varepsilon)=\{
(t,x',x^{d+1}):(t,x')\in \bar\cM_{T} ,
x^{d+1}=0,\pm\varepsilon,\pm2\varepsilon,...\}.
$$
For
$$
\tilde{Q}:=\{(t,x',x^{d+1}):(t,x')\in Q,
x^{d+1}=0,\pm\varepsilon,...,\pm m\varepsilon\},
$$
we have
$$
\tilde{Q}^{o}_{\varepsilon}=Q^{o}_{0}\times\{0,\pm\varepsilon,
...,\pm(m-1)\varepsilon\},
$$
$$
\partial_{\varepsilon}\tilde{Q}
=(\partial_{0}Q\times\{0,\pm\varepsilon,
...,\pm m\varepsilon\})\cup
(Q^{o}_{0}\times\{m\varepsilon,-m\varepsilon\}).
$$

Next, define
$$
\tilde{a}^{\alpha}_{k}(t,x',x^{d+1})
=\begin{cases} a^{\alpha}_{k}(t,x')\quad\text{if}&x^{d+1}>0,\\
\hat{a}^{\alpha}_{k}(t,x')\quad\text{if}&x^{d+1}\leq0,
\end{cases}
$$
and similarly introduce $\tilde{b}^{\alpha}_{k}$ and
$\tilde{c}^{\alpha}$. Let
$$
\tilde{f}^{\alpha}(t,x',x^{d+1})
=\begin{cases} f^{\alpha}(t,x')[1-(x^{d+1}-\varepsilon)
/(\varepsilon m)]
&\text{if}\quad x^{d+1}>0,\\
\hat{f}^{\alpha}(t,x')[1+x^{d+1}/(\varepsilon m)]
&\text{if}\quad x^{d+1}\leq0,
\end{cases}
$$
and similarly define $\tilde{u}(t,x',x^{d+1})$.

Next, we check that   Theorem \ref{theorem 8.19.1}
is applicable
to $\tilde{Q},\tilde{u},\tilde{a},\tilde{b},\tilde{c}$,
and $\tilde{f}$.
Obviously, $\tilde{u}$ in $\tilde{Q}
\cap[0,T)\times\bR^{d+1}$ satisfies 
equation \eqref{8.19.1}
constructed on the basis of $\tilde{a},\tilde{b},\tilde{c}$,
and $\tilde{f}$. In Assumption \ref{assumption 8.25.2}
inequalities (\ref{8.24.03}) and (\ref{8.25.3})
for $\delta_{h,\ell_{k}}$ hold by assumption.
To check them for $\delta_{\varepsilon,\pm l}$, observe that
in $\tilde{Q}^{o}_{\varepsilon}$
$$
\delta_{\varepsilon,l}(\tilde{a}^{\alpha}_{k},
\tilde{b}^{\alpha}_{k},\tilde{c}^{\alpha},\tilde{f}^{\alpha})
(t,x)=\begin{cases}(0,0,0,-f^{\alpha}(t,x')/ m )&\quad\text{if}
\quad x^{d+1}>0,\\
(0,0,0,\hat f^{\alpha}(t,x')/  m)&\quad\text{if}
\quad x^{d+1}<0,
\end{cases}
$$
and
$$
\delta_{\varepsilon,l}(\tilde{a}^{\alpha}_{k},
\tilde{b}^{\alpha}_{k},\tilde{c}^{\alpha},\tilde{f}^{\alpha})
(t,x',0)=\varepsilon^{-1}
(a^{\alpha}_{k}-\hat{a}^{\alpha}_{k},
b^{\alpha}_{k}
-\hat{b}^{\alpha}_{k},c^{\alpha}-\hat{c}^{\alpha},
f^{\alpha}-\hat{f}^{\alpha})(t,x'),
$$
where by virtue of (\ref{8.26.2})
$$
\varepsilon^{-1}
|a^{\alpha}_{k}(t,x')-\hat{a}^{\alpha}_{k}(t,x')|\leq
K\sqrt{\tilde{a}^{\alpha}_{k}(t,x',0)}+Kh.
$$
Using the above formulas along with (\ref{8.26.1})
and the inequality $\varepsilon m\geq1$ we conclude that
in our situation
 (\ref{8.24.03}) and (\ref{8.25.3})  hold  
 for $\delta_{\varepsilon, l}$. The same is true
for
$\delta_{\varepsilon,-l}=-T_{\varepsilon,-l}\delta
_{\varepsilon,l}$ .

Now by Theorem \ref{theorem 8.19.1} we obtain that
for $(t,x')\in Q$
$$
\varepsilon^{-1}|u(t,x')-\hat{u}(t,x')|
=|\delta_{\varepsilon,l}\tilde{u}(t,x',0)|
 \leq 
N e^{c_{0}(T+\tau)}\big(1+|u|_{0,Q}+|\hat u|_{0,Q}
$$
\begin{equation}
                                                 \label{8.26.3}
+\max_{\partial_{0}Q}(\max_{k}|\delta_{h ,\ell_{k}}u|
+\max_{k}|\delta_{h ,\ell_{k}}\hat{u}|
+\varepsilon^{-1}|u-\hat{u}|)+I_{m}\big),
\end{equation}
where 
$$
I_{m}:=\max_{Q^{o}_{0}\times\{m\varepsilon,-m\varepsilon\}}
(\max_{k}|\delta_{h ,\ell_{k}}\tilde{u}|
+|\delta_{\varepsilon,l}\tilde{u}|+
|\delta_{\varepsilon,-l}\tilde{u}|).
$$
Since on $Q^{o}_{0}\times\{r\varepsilon:r=m,m\pm1\}$
$$
|\tilde{u}(t,x)|\leq N(1/m+|1-|x^{d+1}|/(\varepsilon m)|)
\leq N/m,
$$
where $N$ is independent of $m$, we have $I_{m}\to0$
as $m\to\infty$
and by letting $m\to\infty$ in \eqref{8.26.3}, we arrive
at \eqref{8.25.1}. The theorem is proved.

We also need a version of Theorem \ref{theorem 8.24.1}
in the case that $Q=\bar\cM_{T}$.
In the following theorem we 
abandon Assumption \ref{assumption 8.25.3} and
go back to our basic assumptions.

\begin{theorem}
                                       \label{theorem 9.11.3}
Let $\hat{\sigma}^{\alpha}_{k},
\hat{b}^{\alpha}_{k},\hat{c}^{\alpha},\hat{\lambda},
\hat{f}^{\alpha}$ satisfy the assumptions
 in Section~\ref{section 10.4.1} and $\hat{\lambda}=\lambda$.
 Let
$u$  be a function on $\bar\cM_{T}$
satisfying  \eqref{8.19.1}  in $\cM_{T}$ and
let   $\hat{u}$ be a function  on $\bar\cM_{T}$  satisfying equation
  \eqref{8.19.1}  in $\cM_{T}$ with 
$\hat{a}^{\alpha}_{k},
\hat b^{\alpha}_{k},\hat c^{\alpha}, 
\hat f^{\alpha}$   in place of
$a^{\alpha}_{k},
b^{\alpha}_{k},c^{\alpha},  
f^{\alpha}$, respectively. Assume that
$u$ and $\hat{u}$ are bounded
on $\bar\cM_{T}$ and  
$$
|u(T,\cdot)|,|\hat{u}(T,\cdot)|\leq 
K.$$
Introduce
$$
\varepsilon=
\sup_{\cM_{T},A,k}
\big(|\sigma^{\alpha}_{k}-\hat{\sigma}^{\alpha}_{k}|+
|b^{\alpha}_{k}-\hat{b}^{\alpha}_{k}|
+|c^{\alpha} -\hat{c}^{\alpha} |
+|f^{\alpha} -\hat{f}^{\alpha} |\big).
$$

 Then there is a constant $N$ depending only
on $K $  and $d_{1}$ such that   if for a number
$c_{0}\geq0$ equation \eqref{9.11.11} holds,
then
\begin{equation}
                                                  \label{8.25.01}
| u-\hat{u}| \leq 
N\varepsilon e^{c_{0}(T+\tau)}I
\end{equation}
on $\bar\cM_{T}$, where
$$
I=\sup_{(T,x)\in\bar\cM_{T}}\big(1 
+(\max_{k}|\delta_{h ,\ell_{k}}u|
+\max_{k}|\delta_{h ,\ell_{k}}\hat{u}|
+\varepsilon^{-1}|u-\hat{u}|)(T,x)\big).
$$

\end{theorem}

Proof.
First we show that we may assume that
$\varepsilon\in(0, h]$. To this end for $\theta\in[0,1]$
introduce $u^{\theta}$ as the unique
bounded solution
of
$$
\delta^{T}_{\tau}u
+\sup_{\alpha\in A}[a^{\theta\alpha}_{k}
\Delta_{h,\ell_{k}}u
+ b^{\theta\alpha}_{k} 
\delta_{h,\ell_{k}}u
+c^{\theta\alpha}  u
+f^{\theta\alpha} ]=0
$$
in $\cM_{T}$
 and $u=(1-\theta)u+\theta\hat{u}$ in $\{(T,x)\in\cM\}$, where
$$
[\sigma^{\theta\alpha}_{k}, b^{\theta\alpha}_{k},
c^{\theta\alpha},f^{\theta\alpha} ]=
(1-\theta)[\sigma^{\alpha}_{k},b^{\alpha}_{k},c^{\alpha} ,
f^{\alpha} ]+\theta[\hat
\sigma^{\alpha}_{k},\hat b^{\alpha}_{k},\hat c^{\alpha} ,
\hat f^{\alpha}],
$$
$$
a^{\theta\alpha}_{k}=
(1/2)|\sigma^{\theta\alpha}_{k}|^{2}.
$$

 Obviously, $u^{0}=u$ and $u^{1}=\hat{u}$. Also
notice that for any $\theta_{1},\theta_{2}\in[0,1]$
$$
|\sigma^{\theta_{1}\alpha}_{k}-\sigma^{\theta_{2}\alpha}_{k}|+
|b^{\theta_{1}\alpha}_{k}-b^{\theta_{2}\alpha}_{k}|
+|c^{\theta_{1}\alpha}-c^{\theta_{2}\alpha}|
+|c^{\theta_{1}\alpha}-c^{\theta_{2}\alpha}|
\leq  |\theta_{1}-\theta_{2}|\varepsilon.
$$

Therefore, if the present theorem holds true for
$\varepsilon\in(0, h]$, then
for any $\varepsilon>0$ as long as
$|\theta_{1}-\theta_{2}|\varepsilon\leq  h$ we have
$$
| u^{\theta_{1}}-u^{\theta_{2}}| \leq 
N_{1}|\theta_{1}-\theta_{2}|
\varepsilon e^{c_{0}T}I(\theta_{1},\theta_{2}),
$$
where
$$
I(\theta_{1},\theta_{2})=
\sup_{(T,x)\in\bar\cM_{T}}\big(1 
+(\max_{k}|\delta_{h ,\ell_{k}} u^{\theta_{1}}|
+\max_{k}|\delta_{h ,\ell_{k}} u^{\theta_{2}}|
$$
$$
+|\theta_{1}-\theta_{2}|^{-1}
\varepsilon^{-1}| u^{\theta_{1}}- u^{\theta_{2}}|)(T,x)\big).
$$
Obviously, $I(\theta_{1},\theta_{2})\leq 4I$, so that
$$
| u^{\theta_{1}}-u^{\theta_{2}}| \leq 
4N_{1}|\theta_{1}-\theta_{2}|
\varepsilon e^{c_{0}T}I .
$$
By dividing the interval $(0,1)$
into pieces of appropriate length and adding up these estimates
we come to \eqref{8.25.01} with the constant $N$ which is
$4$ times larger than the one which suits $\varepsilon
\leq  h$.

Thus indeed the only important case is the one with
$\varepsilon\in(0, h]$. In this case,
actually, the theorem is a simple consequence
of Theorem~\ref{theorem 8.24.1}, Corollaries
\ref{corollary 9.9.2} and
\ref{corollary 9.9.1} and Lemma \ref{lemma 9.9.3}.
Indeed, by Lemma \ref{lemma 9.9.3} we can approximate
both $u$ and $\hat{u}$ with solutions
such that $f$ and $\hat{f}$ have compact support
as well as the restriction of approximating functions
to $\{t=T\}$. For approximating functions
we get the result as in the proof of 
Theorem \ref{theorem 8.24.1} by expanding
finite sets $Q$ and using that the contribution
coming from the distant
boundary becomes negligible due
to Corollary
\ref{corollary 9.9.1}. 
 We also get rid
of terms $|u|_{0,Q}$ and $|\hat{u}|_{0,Q}$
on the basis of Corollary
\ref{corollary 9.9.2}. However, to use 
Theorem \ref{theorem 8.24.1} we also have
to notice that due to the assumption that
$\varepsilon\leq h$ we have
$$
|a^{ \alpha}_{k}-\hat a^{ \alpha}_{k}|
\leq(|\sigma^{ \alpha}_{k}|\wedge |\hat\sigma^{ \alpha}_{k}|) 
|\sigma^{ \alpha}_{k}- \hat\sigma^{ \alpha}_{k}|
+|\sigma^{ \alpha}_{k}
- \hat\sigma^{ \alpha}_{k}|^{2}
\leq 2 \varepsilon
\sqrt{a^{ \alpha}_{k}\wedge \hat a^{ \alpha}_{k}}
+  \varepsilon^{2}
$$
and $\varepsilon^{2}\leq\varepsilon h$.
The theorem is proved.

\mysection{H\"older continuity of $v$ and $v_{\tau,h}$
in $t$}
                                       \label{section 10.2.5}

We will be using the method of ``shaking'' the coefficients
introduced in \cite{Kr99} and \cite{Kr00}.
Take a nonempty set 
$$
S\subset B_{1}=\{x\in\bR^{d}:|x|<1\}
$$
and for $\varepsilon\in\bR^{d}$ introduce
$v_{\tau,h}^{\varepsilon,S}$ as the unique
solution of equation
$$
\delta^{T}_{\tau}u+\sup_{
(\alpha,y)\in A\times S}[
L^{\alpha}_{h}(t,x+\varepsilon y)u(t,x)
+f^{\alpha}(t,x+\varepsilon y)]=0
$$
in $H_{T}$ with terminal condition
\begin{equation}
                                             \label{9.14.2}
u(T,x)=\sup_{y\in S}g(x+\varepsilon y)\quad\text{on}\quad
 \bR^{d}.
\end{equation}
Also let $v ^{\varepsilon,S}$
be a probabilistic
solution of
$$
\frac{\partial}{\partial t} u(t,x)+
\sup_{
(\alpha,y)\in A\times S}[
L^{\alpha} (t,x+\varepsilon y)u(t,x)
+f^{\alpha}(t,x+\varepsilon y)]
=0
$$
in $H_{T}$ with terminal condition
\eqref{9.14.2}.
  Observe that if $S$ is a singleton
$\{y\}$, then by uniqueness
$$
 v_{\tau,h}^{\varepsilon,S}(t,x)=v_{\tau,h}(t,x+\varepsilon y),
\quad v ^{\varepsilon,S}(t,x)=v (t,x+\varepsilon y).
$$

\begin{lemma}
                                    \label{lemma 9.11.1}
There is a constant $N$ depending only on
$K$ and $d_{1}$ such that
if for a number
$c_{0}\geq0$ equation \eqref{9.11.11} holds, 
then for all $\varepsilon
\in\bR$ 
\begin{equation}
                                                    \label{9.11.6}
|v_{\tau,h}^{\varepsilon,S}-v_{\tau,h}|
\leq Ne^{c_{0}(T+\tau)}|\varepsilon|
\quad\text{on}\quad\bar H_{T}, 
\end{equation}
\begin{equation}
                                                    \label{9.11.7}
|v ^{\varepsilon,S}-v |
\leq Ne^{ (N -\lambda)_{+}T}|\varepsilon|\quad\text{on}\quad\bar H_{T}.
\end{equation}
In particular, (take $S=\{(y-x)/|y-x|\}$, $\varepsilon=|y-x|$)
$$
|v_{\tau,h}(t,y)-v_{\tau,h}(t,x)|\leq N
e^{c_{0}(T+\tau)}|y-x|,
\quad (t,y),(t,x)\in\bar H_{T}, 
$$
\begin{equation}
                                                    \label{9.11.9}
|v (t,y)-v (t,x)|\leq Ne^{ (N -\lambda)_{+}T}|y-x|,\quad
(t,y),(t,x)\in \bar H_{T}.
\end{equation}
\end{lemma}

Proof. While proving \eqref{9.11.6} 
we may concentrate on $\bar\cM_{T}$. Then
it suffices to use 
 Theorem \ref{theorem 9.11.3},
where we   take $A\times S$, $(\sigma,b,c,f)(t,x  )$
and $(\sigma,b,c,f)(t,x+\varepsilon y)$ in place of $A$,
$(\sigma,b,c,f) $ and $(\hat{\sigma},\hat{b},\hat{c},\hat{f}) $,
respectively. We also use that the difference of sups
is less than the sup of differences
while estimating the boundary terms.

Estimate \eqref{9.11.9} is a particular case
of Theorem 4.1.1 of \cite{Kr77} and \eqref{9.11.7}
is, actually, a particular case of \eqref{9.11.9}
since one can view $\varepsilon$ as just another
coordinate of the space variable. The lemma is proved.

For $\Lambda\subset(-1,0)$ introduce
$v_{\tau,h}^{\varepsilon,\Lambda,S}$ as the unique
bounded
solution of equation
$$
\delta^{T}_{\tau}u(t,x)+\sup_{
(\alpha,r,y)\in A\times\Lambda\times S}[
L^{\alpha}_{h}(t+\varepsilon^{2}r,x+\varepsilon y)u(t,x)
$$
\begin{equation}
                                             \label{9.11.02}
+f^{\alpha}(t+\varepsilon^{2}r,x+\varepsilon y)]=0
\end{equation}
in $ H_{T}$ with terminal condition \eqref{9.11.4}.
Also let $v ^{\varepsilon,\Lambda,S}$
be a probabilistic
solution of
$$
\frac{\partial}{\partial t} u(t,x)+
\sup_{
(\alpha,r,y)\in A\times\Lambda\times S}[
L^{\alpha} (t+\varepsilon^{2}r,x+\varepsilon y)u(t,x)
$$
$$
+f^{\alpha}(t+\varepsilon^{2}r,x+\varepsilon y)]
=0
$$
in $  H_{T}$ with terminal condition \eqref{9.11.4}.

\begin{lemma}
                                    \label{lemma 9.13.3}
There is a constant $N$ depending only on
$K$  and $d_{1}$ such that
if for a number
$c_{0}\geq0$ equation \eqref{9.11.11} holds
and
assumption (H) of Theorem \ref{theorem 8.26.1} is satisfied,
 then   for all $\varepsilon
\in\bR$ 
\begin{equation}
                                                    \label{9.13.2}
|v_{\tau,h}^{\varepsilon,\Lambda,S}-v_{\tau,h}|
\leq Ne^{c_{0}(T+\tau)}|\varepsilon|
\quad\text{on}\quad[0,T]\times
\bR^{d}, 
\end{equation}
If, additionally,  $\tau,h\leq K$, then
\begin{equation}
                                                    \label{9.13.6}
|v^{\varepsilon,\Lambda,S}_{\tau,h}(t,x)
-v^{\varepsilon,\Lambda,S}_{\tau,h}(s,y)|\leq N 
 e^{c_{0}(T+\tau)}(|t-s|^{1/2}+|y-x| ),
\end{equation}
\begin{equation}
                                                    \label{9.13.4}
|v_{\tau,h}(t,x)-v_{\tau,h}(s,y)|\leq  N 
e^{c_{0}(T+\tau)}( |t-s|^{1/2}+|y-x| )
\end{equation}
for all $(t,x),(s,y)\in\bar H_{T}$ with $|t-s|\leq1$,
\begin{equation}
                                                    \label{9.13.3}
|v ^{\varepsilon,\Lambda,S}-v |
\leq  N e^{ (N -\lambda)_{+}T}|\varepsilon|\quad\text{on}\quad[0,T]\times
\bR^{d},
\end{equation}
\begin{equation}
                                                    \label{9.13.7}
|v^{\varepsilon,\Lambda,S}(t,x)
-v^{\varepsilon,\Lambda,S}(s,y)|\leq  N 
e^{ (N -\lambda)_{+}T}(|t-s|^{1/2}+|y-x|) ,
\end{equation}
\begin{equation}
                                                    \label{9.13.5}
|v (t,x)-v (s,y)|\leq  N 
e^{ (N -\lambda)_{+}T}(|t-s|^{1/2}+|y-x|) 
\end{equation}
for all $(t,x),(s,y)\in\bar H_{T}$ with $|t-s|\leq1$.

\end{lemma}

Proof. Estimates \eqref{9.13.3} and
 \eqref{9.13.7} are proved in Corollary 3.2 of \cite{Kr99}.
Estimate \eqref{9.13.5} is obtained from
 \eqref{9.13.7} by setting $\varepsilon=0$.

The proof of \eqref{9.13.2} follows that of 
\eqref{9.11.6}  and is left to the reader.

On the one hand, estimate \eqref{9.13.6} for $\varepsilon=0$
implies \eqref{9.13.4} and, on the other hand,
\eqref{9.11.02} is a particular case of
\eqref{8.19.1}, and therefore 
 \eqref{9.13.6} is a particular case of
\eqref{9.13.4}. Hence to finish proving the lemma
it only remains to prove \eqref{9.13.4}.

Since, as is stated in Lemma \ref{lemma 9.11.1},
 $v_{\tau,h}$ is Lipschitz continuous in $x$,
we   only remains to  prove that 
$$
I(t,s,x):=|v_{\tau,h}(t,x)-v_{\tau,h}(s,x)|
\leq N|t-s|^{1/2}.
$$

In addition, if $0\leq t\leq s\leq T$,
$s-t\leq1$, and $s-t =n\tau+\gamma$, 
where
$n=0,1,...$, $\gamma\in[0,\tau)$,
 then by  Lemma \ref{lemma 9.11.1} and 
Corollary \ref{corollary 9.9.3}
$$
I(t,s,x)\leq 
|v_{\tau,h}(t,x)-v_{\tau,h}(t+n\tau,x)|+
|v_{\tau,h}(t+ n\tau,x)-v_{\tau,h}(s,x)|
$$
$$
\leq
Ne^{c_{0}(T+\tau)}
|t-s|^{1/2}+|v_{\tau,h}(t+n\tau,x)-v_{\tau,h}(s,x)|.
$$
Thus, it suffices to estimate $I(t,s ,x )$
for $s = t +\gamma$ with $\gamma\in(0,\tau)$.
By shifting the origin we reduce the problem
to showing that
\begin{equation}
                                              \label{9.14.4}
I(0,\gamma,0)\leq Ne^{c_{0}(T+\tau)}\gamma^{1/2}.
\end{equation}

Introduce $S=\tau[T/\tau]$ and
first, additionally assume that $S\geq\tau$.
In that case,  set $u=v_{\tau,h}$,
$\hat u(r,y)=v_{\tau,h}((r+\gamma)\wedge T,y)$, and
$$
[\hat
\sigma^{\alpha}_{k},\hat b^{\alpha}_{k},\hat c^{\alpha} ,
\hat f^{\alpha}](r,y)
=[\sigma^{ \alpha}_{k}, b^{ \alpha}_{k},
c^{ \alpha},f^{ \alpha}](r+\gamma,y).
$$
Notice that  for $(r,y)\in\cM_{S}$
we have $ r+\gamma<S\leq T$,
$$
\tau_{S}(r)=\tau,\quad r+\tau_{S}(r)=r+\tau,\quad
(r+\tau+\gamma)\wedge T=r+\gamma+\tau_{T}(r+\gamma)
$$
$$
\hat{u}(r+\tau_{S}(r),y)-\hat{u}(r,y)
=v_{\tau,h}((r+\tau+\gamma)\wedge T,y)-v_{\tau,h}(r+\gamma,y).
$$
It follows that relative to $\bar\cM_{S}$
the function $\hat u$ in $\cM_{S}$
 satisfies equation \eqref{8.19.1}
constructed from
$\hat
\sigma^{\alpha}_{k},\hat b^{\alpha}_{k},\hat c^{\alpha} ,
\hat f^{\alpha}$. 
By observing that the parameter $\varepsilon$
in Theorem~\ref{theorem 9.11.3} is less
than $K\gamma^{1/2}$ owing to assumption (H)
of Theorem~\ref{theorem 8.26.1} and using again that
$v_{\tau,h}$ is Lipschitz continuous in $x$
we obtain from Theorem \ref{theorem 9.11.3} that
$$
I(0,\gamma,0)=|v_{\tau,h}(0,0)-v_{\tau,h}(\gamma,0)|
=|u(0,0)-\hat{u}(0,0)|
$$
$$
\leq
Ne^{c_{0}(T+\tau)}\gamma^{1/2}
+\sup_{(S,y)\in\bar\cM_{S}}|u(S,y)-\hat{u}(S,y)|
$$
$$
=Ne^{c_{0}(T+\tau)}\gamma^{1/2}+\sup_{y}
|v_{\tau,h}(S,y)-v_{\tau,h}((S+\gamma)\wedge T,y)|.
$$
Thus, after one more shift of the origin,
bringing $S$ to zero,
we reduce the problem of estimating $I(0,\gamma,0)$
to the situation
when $T<\tau$, so that
$t=0$, $\tau (t)=T-t$, and $t+\tau(t)=T$ on 
$$
 \cM_{T}=\bar\cM_{T}\cap\{t=0\}.
$$

 Then the function $\tilde{u}$,
introduced  on $ \bar\cM_{T}$ by
$$
\tilde{u}(0,x)=v_{\tau,h}(\gamma,x),\quad
\tilde{u}(T,x)=g(x),
$$
on $\cM_{T}$ satisfies equation \eqref{8.19.1}
corresponding to
$\hat
\sigma^{\alpha}_{k},\hat b^{\alpha}_{k},\hat c^{\alpha} ,
\hat f^{\alpha}$. By Theorem \ref{theorem 9.11.3}
we conclude that
$$
I(0,\gamma,0)= |v_{\tau,h}(\gamma,0)
-v_{\gamma,h}(0,0)|= |\tilde{u}(0,0)-u(0,0)|\leq N\gamma^{1/2}.
$$
Estimate \eqref{9.14.4} 
  and the lemma are proved.

\mysection{Proof of Theorems \protect\ref{theorem 8.26.1},
\protect\ref{theorem 8.27.1}, and \protect\ref{theorem 9.19.1}}
                                        \label{section 10.2.1}

{\bf Proof of Theorem  \ref{theorem 8.26.1}}.
We start with proving  \eqref{9.18.2}
with $N$ which may depend on $T$. Observe that if 
$$
T\leq2\varepsilon^{2},\quad \varepsilon:= (\tau+h^{2})^{1/4},
$$ 
then we have nothing to prove
since then by \eqref{9.13.5}
and \eqref{9.13.4}
$$
\sup_{\bar H_{T}}|v_{\tau,h}-v|
\leq\sup_{\bar H_{T}}|v_{\tau,h}-g|
+\sup_{\bar H_{T}}|v_{\tau,h}-g|
\leq NT^{1/2}
$$
$$
\leq N(\tau+h^{2})^{1/4}\leq N(\tau^{1/4}+h^{1/2}).
$$

Therefore in the rest of the proof
without losing generality we assume that $T>2\varepsilon^{2}$.
By Corollary \ref{corollary 9.9.2} we have $|v|$ and
$|v_{\tau,h}|$ under control and therefore
  we may assume that $ h\leq1$
and $\tau$ is so small that there is a $c_{0}
=c_{0}(K,d_{1})$ such that even with $\lambda=0$
it satisfies condition
\eqref{9.11.11} imposed in Lemma 
\ref{lemma 9.13.3}.

First we prove that
\begin{equation}
                                              \label{9.15.1}
v\leq v_{\tau,h}+N(\tau^{1/4}+h^{1/2})\quad\text{on}
\quad\bar H_{T}.
\end{equation}

We take $\Lambda=(-1,0)$ and $S=B_{1}$ 
and set
$$
v^{\varepsilon}_{\tau,h}=v^{\varepsilon,\Lambda,S}_{\tau,h},
$$
where the latter function is introduced before
Lemma \ref{lemma 9.13.3}. Then for any $\alpha\in A$,
$r\in(-1,0)$, and  $|y|<1$ 
\begin{equation}
                                              \label{9.15.2}
\delta_{\tau}v^{\varepsilon}_{\tau,h}(t-\varepsilon^{2}r,
x-\varepsilon y)+L^{\alpha}_{h}(t,x)
v^{\varepsilon}_{\tau,h}(t-\varepsilon^{2}r,
x-\varepsilon y)+f^{\alpha}(t,x)\leq0
\end{equation}
 provided that
$$
(t,x)\in\bar H_{T-2\varepsilon^{2}}
\subset
\bar H_{T-\tau-\varepsilon^{2}}.
$$

Next take
  a nonnegative function $\zeta\in C_{0}^{\infty}(\bR^{d+1})$
with support in $(-1,0)\times B_{1}$
and unit integral. For any function
$u$ for which it makes sense we set
$$
u^{(\varepsilon)}(t,x)=\varepsilon^{-d-2}\int
_{\bR^{d+1}}u(s,y)\zeta((t-s)/\varepsilon^{2},(x-y)
/\varepsilon)\,dsdy.
$$
By multiplying \eqref{9.15.2} by $\zeta$ and integrating
we get that for any $\alpha\in A$ on 
$\bar H_{T-2\varepsilon^{2}}$ it holds that
$$
\delta_{\tau}v^{\varepsilon(\varepsilon)}_{\tau,h}
+L^{\alpha}_{h}v^{\varepsilon(\varepsilon)}_{\tau,h}
+f^{\alpha}\leq0.
$$
From here by Taylor's formula (cf.~\eqref{9.9.2}) we infer
$$
\frac{\partial}{\partial t}
v^{\varepsilon(\varepsilon)}_{\tau,h}
+L^{\alpha}v^{\varepsilon(\varepsilon)}_{\tau,h}
+f^{\alpha}\leq N\big(\tau|D^{2}_{t}
v^{\varepsilon(\varepsilon)}_{\tau,h}|
_{0,\bar H_{T-2\varepsilon^{2}}}
$$
$$
+h^{2}|D^{4}_{x}
v^{\varepsilon(\varepsilon)}_{\tau,h}|_{0,
\bar H_{T-2\varepsilon^{2}}}
+h|D^{2}_{x}
v^{\varepsilon(\varepsilon)}_{\tau,h}|_{0, 
\bar H_{T-2\varepsilon^{2}}}\big)=:I
$$
in $\bar H_{T-2\varepsilon^{2}}$.
It follows that 
\begin{equation}
                                              \label{10.20.1}
v^{\varepsilon(\varepsilon)}_{\tau,h}+
(T-2\varepsilon^{2}-t)I
\end{equation}
is a supersolution of \eqref{8.19.3} in 
$\bar H_{T-2\varepsilon^{2}}$ and either by It\^o's
formula or by   properties of viscosity
solutions we have in $\bar H_{T-2\varepsilon^{2}}$
that
\begin{equation}
                                              \label{9.15.3}
v\leq v^{\varepsilon(\varepsilon)}_{\tau,h}+(T
-2\varepsilon^{2}-t)I
+\sup_{\{T-2\varepsilon^{2}\}\times\bR^{d}}
|v-v^{\varepsilon(\varepsilon)}_{\tau,h}|.
\end{equation}

Now use the fact that  owing to \eqref{9.13.7}
and well-known properties of convolutions
we have in $\bar H_{T-2\varepsilon^{2}}$ that
$$
|v^{\varepsilon(\varepsilon)}_{\tau,h}-
v^{\varepsilon }_{\tau,h}|\leq N\varepsilon
$$
with $N$ depending only on $K$, $T$, $d$, and $d_{1}$
and
for any $n=1,2,...$
$$
|D^{n}_{t}v^{\varepsilon(\varepsilon)}_{\tau,h}|
_{0,\bar H_{T-2\varepsilon^{2}}}
+|D^{2n}_{x}v^{\varepsilon(\varepsilon)}_{\tau,h}|
_{0,\bar H_{T-2\varepsilon^{2}}}
\leq N/\varepsilon^{2n-1},
$$
where $N$ depends only on $n$, $K$, $T$, $d$, and $d_{1}$.
Also, notice that 
$$
|v(T-2\varepsilon^{2},x)-
v^{\varepsilon }_{\tau,h}(T-2\varepsilon^{2},x)|,
$$
that appears from the last term in \eqref{9.15.3},
is estimated through $N\varepsilon$
  in the beginning of the proof.
 Then we conclude
$$
 v\leq v_{\tau,h}+N[\varepsilon +
(\tau+h^{2})/\varepsilon^{3}+ h/\varepsilon]
$$
 in $\bar H_{T-2\varepsilon^{2}}$. Actually,
the same estimate holds in
$\bar H_{T}$ due to the argument in
the beginning of the proof. Finally by observing that
$$
\varepsilon +
(\tau+h^{2})/\varepsilon^{3}+ h/\varepsilon \leq
 \varepsilon+(\tau+h^{2})/\varepsilon^{3}
+(\tau+h^{2})^{1/2}/\varepsilon 
$$
and 
recalling that $\varepsilon=(\tau+h^{2})^{1/4}$ we come to
\eqref{9.15.1}.

It remains to prove that
\begin{equation}
                                              \label{9.15.4}
v_{\tau,h}-v\leq N(\tau^{1/4}+h^{1/2}).
\end{equation}

Similarly to what was done with the discrete
approximation above, on the basis
of functions $v ^{\varepsilon,\Lambda,S}$
in the proof of Theorem 2.1 of \cite{Kr00}  
an   infinitely
differentiable  function $u$
on $\bar H_{T}$ is constructed such that
$$
\frac{\partial}{\partial t}u
+\sup_{\alpha\in A}[L^{\alpha}u+f^{\alpha}]\leq 0,
\quad |u-v|\leq N\varepsilon\quad\text{on}\quad
\bar H_{T},
$$
with $N$ depending only on $K$, $T$, $d$, and $d_{1}$
and for any $n=1,2,...$
$$
|D^{n}_{t}u|_{0,\bar H_{T}}
+|D^{2n}_{x}u|_{0,\bar H_{T}}
\leq N/\varepsilon^{2n-1},
$$
where $N$ depends only on $n$, $K$, $T$, $d$, and $d_{1}$.
As above, it follows by Taylor's formula   that
on $\bar H_{T-\tau}$ (where $\tau_{T}(t)=\tau$)
we have
$$
\delta^{T}_{\tau}u
+\sup_{\alpha\in A}[L_{h}^{\alpha}u+f^{\alpha}]\leq
N(\tau+h^{2})/\varepsilon^{3}+Nh/\varepsilon.
$$
   Upon taking
$$
u_{1}=v_{\tau,h},\quad u_{2}=u+\sup_{H_{T}\setminus H_{T-\tau}}
(v_{\tau,h}-u)_{+},\quad C=
N(\tau+h^{2})/\varepsilon^{3}+Nh/\varepsilon
$$
in Lemma \ref{lemma 9.9.2}, we obtain
$$
v_{\tau,h}\leq u+\sup_{H_{T}\setminus H_{T-\tau}}
(v_{\tau,h}-u)_{+}
+N(\tau+h^{2})/\varepsilon^{3}+Nh/\varepsilon.
$$
Here $u\leq v+N\varepsilon$ and, owing to \eqref{9.13.5}
and \eqref{9.13.4}
and the above mentioned properties of $u$,
$$
\sup_{H_{T}\setminus H_{T-\tau}}
(v_{\tau,h}-u)_{+}\leq\sup_{H_{T}\setminus H_{T-\tau}}
|v_{\tau,h}-g|
+\sup_{H_{T}\setminus H_{T-\tau}}
|g-v|
$$
$$
+\sup_{H_{T}\setminus H_{T-\tau}}
|v-u|\leq N(\tau^{1/2}+\varepsilon).
$$

Thus,
$$
v_{\tau,h}\leq v+N[\varepsilon+\tau^{1/2}+
(\tau+h^{2})/\varepsilon^{3}+ h/\varepsilon]
$$
$$
\leq N[\varepsilon+(\tau+h^{2})/\varepsilon^{3}
+(\tau+h^{2})^{1/2}/\varepsilon].
$$
Recalling that $\varepsilon=(\tau+h^{2})^{1/4}$ yields
\eqref{9.15.4} and \eqref{9.18.2}
with $N$ perhaps depending on $T$.

However, if $\lambda$ is large enough,
$c_{0}=0$ satisfies condition
\eqref{9.11.11} imposed in Lemma 
\ref{lemma 9.13.3} and for any $\lambda>0$,
the functions $v$ and $v_{\tau,h}$
are bounded by a constant depending only
on $K$ and $\lambda$ owing to
Corollary \ref{corollary 9.9.2}. 
In that case also 
  the estiamtes in Lemma \ref{lemma 9.13.3}
are independent of $T$. Furthermore,
we can replace $T-2\varepsilon^{2}-t$ in
\eqref{10.20.1} with the constant $N$ from
\eqref{9.11.11}.
This allows us
to check that in the above proof the constants
are actually independent of $T$ if
$\lambda\geq N=N(K,d_{1})$.

The theorem is proved.

{\bf Proof of Theorem  \ref{theorem 8.27.1}}.
Take $g\equiv0$ and denote the functions
$v$ and $v_{\tau,h}$ from Theorem \ref{theorem 8.26.1}
by $v^{T}$ and $v^{T}_{\tau,h}$.
 Obviously, it suffices to prove that for all $(t,x)$
\begin{equation}
                                            \label{9.19.2}
 \tilde v(x)=\lim_{T\to\infty}v^{T}(t,x),\quad
\tilde v_{h}(x)=\lim_{T\to\infty}v^{T}_{\tau,h}(t,x),
\end{equation}
whenever $\lambda>0$ and $\tau $ is small enough.

The first relation in \eqref{9.19.2} is well
known  (see, for instance, \cite{FS} or \cite{Kr77}).
To prove the second, it suffices to prove
that for any sequence $T_{n}\to\infty$
such that $v^{T_{n}}_{\tau,h}(t,x)$
converges at all points of $\cM_{\infty}$,
the limit is independent of $t$ and satisfies
\eqref{8.27.2} on
the grid
$$
G=\{i_{i}h\ell_{1}+...+i_{d_{1}}h\ell_{d_{1}}:
i_{k}=0,\pm1,...,k=1,...,d_{1}\}.
$$
Given the former, the latter is obvious. Also notice
that the translation    $t\to t+\tau$ 
brings any solution
of \eqref{8.19.1} on $\cM_{\infty}$   again to  a solution.
Therefore,
it only remains to prove uniqueness
of bounded solutions of 
\eqref{8.19.1}
on $\cM_{\infty}$.

Observe that if $u_{1}$ and $u_{2}$ are two solutions
of \eqref{8.19.1}
on $\cM_{\infty}$,
then they also solve \eqref{8.19.1} on $\cM_{T}$
for any $T$
with terminal condition $u_{1}$ and $u_{2}$,
respectively. By the comparison result
$$
|u_{1}-u_{2}|\leq e^{-\lambda T/2}
\sup|u_{1}-u_{2}|
$$
if $\tau$ is small enough.
Sending $T\to\infty$ proves the uniqueness
and the theorem.

{\bf Proof of Theorem  \ref{theorem 9.19.1}}.
The unique solvability of \eqref{9.19.1}-\eqref{9.11.4}
in the space of bounded functions
is shown be rewriting the problem as
\begin{equation}
                                                 \label{9.19.4}
u(t,x)=g(x)+\int_{t}^{T}
F(\Delta_{h,\ell_{k}}u(s,x),
\delta_{h,\ell_{k}}u(s,x),u(s,x),s, x)\,ds
\end{equation}
and using, say the method of successive approximations.

Next, since $v_{\tau,h}$ are H\"older continuous
in $(t,x)$, for any
sequence $\tau_{n}\downarrow0$, one can find
a subsequence $\tau_{n'}\downarrow0$ such that
$v_{\tau_{n'},h}(t,x)$ converge at each point of $\bR^{d}$
uniformly in $t\in[0,T]$.
Call $u$ the limit of one of subsequences
and introduce
$$
\kappa_{n'}(t)=i\tau_{n'}\quad\text{for}\quad
i\tau_{n'}\leq t<(i+1)\tau_{n'},\quad i=0,1,...
$$
Then for any smooth $\psi(t)$ vanishing
at $t=T$ and $t=0$
$$
\int_{0}^{T}[\psi 
F(\Delta_{h,\ell_{k}}v_{\tau_{n'},h},
\delta_{h,\ell_{k}}v_{\tau_{n'},h}
 , v_{\tau_{n'},h})](\kappa_{n'}(t),x)\,dt
$$
$$
=\int_{0}^{T} v_{\tau_{n'},h}(\kappa_{n'}(t),x)
\frac{\psi(\kappa_{n'}(t),x)-\psi(\kappa_{n'}(t)-\tau_{n'},x)}
{\tau_{n'}}\,dt.
$$
Since the integrands converge uniformly on $[0,T]$
to their natural limits,
we conclude that $u$ satisfies \eqref{9.19.1}
in the weak sense. This is also a continuous function
and $u(T,x)=g(x)$. It follows that $u$
satisfies \eqref{9.19.4} and by uniqueness
$u=v_{h}$. Now Theorem \ref{theorem 9.19.1}
follows directly
from Theorem~\ref{theorem 8.26.1}. 

\mysection{Concluding remarks}
                                       \label{section 9.19.1}

The methods of this article can also be applied to
equations in cylinders like $Q=[0,T)\times D$,
where $D$ is a domain in $\bR^{d}$. 
It is natural to consider \eqref{8.19.3}
and \eqref{8.19.1} in $Q$ with terminal condition
$u(T,x)=g(x)$ in $D$ and require $v$ and $v_{\tau,h}$
be zero in $[0,T]\times(\bR^{d}\setminus D)$.
If we also   assume that $g=0$ on $\partial D$,
then to carry over our methods we only need
to assume that there is a
sufficiently smooth function $\psi$
such that $\psi>0$ in $D$, $\psi=0$ on $\partial D$,
$|\psi_{x}|\geq1$ on $\partial D$, and $L^{\alpha}\psi
<-1$ in $Q$. The reader who went through our proofs
understands that
the only use of   $\psi$
is in estimating the first order finite-differences
of $v_{\tau,h}$ near the lateral boundary of $Q$
and the gradient of $v$ on the lateral boundary of $Q$.

Elliptic problems and semidiscretization can also
be considered in domains. Although these generalizations
are almost straightforward, some additional work yet
needs to be done and to not overburden the present article
with technicalities we decided to put them in
a subsequent article along with a generalization of
Theorem \ref{theorem 9.19.1} to the case 
when assumption (H) is dropped.

Finally, speaking about equations in domains
it is worth noting that
one can reduce a smooth nonzero lateral condition 
to zero just by subtracting the boundary function
from the solution.

\end{document}